\documentclass[11pt]{article}
\usepackage[utf8]{inputenc}
\usepackage{graphicx}
\usepackage{amsmath, amsthm, amssymb} 
\usepackage{hyperref} 
\usepackage{multicol}
\usepackage{multirow}
\usepackage{caption}
\usepackage{xcolor}
\usepackage{subcaption}
\usepackage[labelfont=bf]{caption}
\usepackage{pdflscape}
\usepackage[noend]{algpseudocode}

\usepackage{algorithm}

\newcommand{\myalgsheader}[0]

\algnewcommand{\IIf}[1]{\State\algorithmicif\ #1\ \algorithmicthen}
\algnewcommand{\EndIIf}{\unskip\ \algorithmicend\ \algorithmicif}
\algnewcommand{\IElse}[1]{\State\algorithmicelse\ #1\ \algorithmicthen}
\algnewcommand{\IfThenElse}[3]{% \IfThenElse{<if>}{<then>}{<else>}
	\State \algorithmicif\ #1\ \algorithmicthen\ #2\ \algorithmicelse\ #3}

\begin{document}
	\title{Mixed Precision Rayleigh Quotient Iteration for Total Least Squares Problems}
	
	\author{Eda Oktay\footnote{Department of Numerical Mathematics, Faculty of Mathematics and Physics, Charles University, \{oktay, carson\}@karlin.mff.cuni.cz. Both authors are supported by the Charles University grant No. SVV-2023-260711, GAUK project No. 202722, the Exascale Computing Project (17-SC-20-SC), a collaborative effort of the U.S. Department of Energy Office of Science and the National Nuclear Security Administration, and by the European Union (ERC, inEXASCALE, 101075632). Views and opinions expressed are those of the
authors only and do not necessarily reflect those of the European
Union or the European Research Council. Neither the European Union nor
the granting authority can be held responsible for them.} \hspace{1pt} and Erin Carson\footnotemark[1]}
	\date{}
	\maketitle

	\paragraph{Abstract}
With the recent emergence of mixed precision hardware, there has been a
renewed interest in its use for solving numerical linear algebra
problems fast and accurately. The solution of total least squares
problems, i.e., solving {$\min_{E,r} \| [E, r]\|_F$} subject to
$(A+E)x=b+${$r$}, arises in numerous applications. Solving this problem requires
finding the smallest singular value and corresponding right singular
vector of $[A,b]$, which is challenging when $A$ is large and sparse. An
efficient algorithm for this case due to Bj\"{o}rck et al. [SIAM J. Matrix Anal. Appl. 22(2), 2000], called RQI-PCGTLS, is based on
Rayleigh quotient iteration coupled with the preconditioned conjugate gradient
method. 

We develop a mixed
precision variant of this algorithm, RQI-PCGTLS-MP, in which up to three different precisions can be used. We assume that the lowest precision is used in the computation of the preconditioner, and give theoretical constraints on how this precision must be chosen to ensure stability. In contrast to standard least squares, for total least squares, the resulting constraint depends not only on the matrix $A$, but also on the right-hand side $b$. We perform a number of numerical
experiments on model total least squares problems used in the literature, which demonstrate that our algorithm can attain the same accuracy as RQI-PCGTLS albeit with a potential convergence delay due to the use of low precision. Performance modeling shows that the mixed precision approach can achieve up to a $4\times$ speedup depending on the size of the matrix and the number of Rayleigh quotient iterations performed. 
	
	\section{Introduction}\label{sec1}

The standard least squares (LS) problem involves solving min$_x\|b-Ax\|_2$, where $A\in\mathbb{R}^{m\times n}$ with $m\geq n$ and rank($A$) = $n$. Using the QR factorization 
\begin{equation*}
A = Q\bar{R} = [Q_1\quad Q_2]\begin{bmatrix}
    R\\
    0
\end{bmatrix},
\end{equation*}
 where $Q\in \mathbb{R}^{m\times m}$ is orthogonal and $R\in \mathbb{R}^{n\times n}$ is upper triangular, the solution to the LS problem is given by $x=U^{-1}Q_1^Tb$ and the residual by $\|b-Ax\|_2=\|Q_2^Tb\|_2.$ The LS problem can also be solved via the normal equations, $A^TAx=A^Tb$, which are equivalent to the augmented system
 \begin{equation*}
     \begin{bmatrix}
         I&A\\
         A^T&0
     \end{bmatrix}\begin{bmatrix}
         r\\
         x
     \end{bmatrix}=\begin{bmatrix}
         b\\
         0
     \end{bmatrix}.
 \end{equation*}

Using mixed precision hardware can provide significant performance benefits in numerical linear algebra computations by reducing both computation and communication costs, which has motivated a flurry of recent work in this area; see, e.g., the surveys \cite{a:21, hm:22}. There are various existing mixed precision approaches to solving least squares problems, namely those based on performing iterative refinement on the augmented system.  
The augmented system formulation is originally due to Bj\"{o}rck \cite{b:67}, who showed that iterative refinement for linear systems could be applied, using extra precision in the residual computation and making use of the QR factorization of $A$ to solve a linear system with the augmented matrix in each iteration. For error bounds and implementation with BLAS operations, see \cite{d:09}. Later, the authors in \cite{chp:20} showed that GMRES-based iterative refinement can also be used, in which preconditioned GMRES is used as the solver within each refinement step. The QR factors of $A$, computed in potentially low precision, are used to construct a left preconditioner for the augmented matrix. This method is called GMRES-LSIR. 

An alternative approach is to use mixed precision iterative refinement on the normal equations \cite{hp:21}.  Since $A^TA$ is symmetric positive definite, the Cholesky factor of this matrix, again potentially computed in lower precision (see also \cite{ytd:15}), can be used to construct a  preconditioner. 

  Although the LS problem is commonly encountered in applications, it may not be a realistic approach in some cases. Since LS problems are based on the standard linear model, $Ax=b+r$, with a random error vector $r$, it is assumed that the only error is in the right-hand side, $b$. However, especially in statistics, the matrix $A$ may also affected by sampling and modeling errors. The errors-in-variables model represents errors in both $A$ and $b$ and is denoted as $(A+E)x=b+r$, where $E$ and $r$ are a random error matrix and vector, respectively.

  The first numerically stable algorithm for solving the errors-in-variables model was introduced in \cite{gv:80}. The algorithm solves the minimization problem $min_{E,r}\|[E,r]\|_F$ subject to $(A+E)x=b+r$. Once a suitable $(E,r)$ pair is found, any solution $x$ to the system $(A+E)x=b+r$ is then considered to be a solution to the total least squares (TLS) problem.

  The exact solution to the TLS problem is given by the singular value decomposition (SVD), $[A,b]=U\Sigma V^T$, where $U$ is $m\times (n+1)$ with orthonormal columns, $\Sigma=$ diag$(\sigma_1,\cdots,\sigma_{n+1})$, and $V=[v_1,\cdots,v_{n+1}]$ is an $(n+1)\times (n+1)$ orthogonal matrix. Assume that the singular values of the matrix $[A,b]$ are $\sigma_1\geq \sigma_2 \geq \cdots \geq \sigma_{n+1} \geq 0$ and the singular values of $A$ are $\sigma'_1\geq \sigma'_2 \geq \cdots \geq \sigma'_n \geq 0.$ If $\sigma_{n+1}=0$, then $[E,r]=0$. Otherwise, {min}$_{\text{rank}(A+E,b+r)<n+1}\|[E,r]\|_F=\sigma_{n+1}$ and the TLS solution is \[x_{TLS}=-\frac{1}{v_{n+1,n+1}}[v_{1,n+1},\cdots,v_{n,n+1}]^T.\] However, computing the SVD explicitly has a high computational cost and is thus not the preferred approach in practice.

  In \cite{gv:80}, the authors show that the approximate condition number for the TLS problem is \[\kappa_{TLS}(A,b) = \dfrac{\sigma'_1}{\sigma'_n-\sigma_{n+1}} = \kappa(A)\dfrac{\sigma'_n}{\sigma'_n-\sigma_{n+1}},\] 
  where $\kappa(A)$ denotes the $2$-norm condition number of $A$. As stated in \cite{bhm:00}, when $1-\sigma_{n+1}/\sigma'_n \ll 1$, $\kappa_{TLS}(A,b)$ can be much greater than $\kappa(A)$. Thus the TLS problem can be very sensitive to perturbations in {the} data.

  To reduce the computational cost of solving the TLS problem, an iterative approach can be used. There are various iterative approaches to solving the TLS problem, such as inverse iteration. However, inverse iteration convergence depends on the ratio $\|\sigma_{n-p}-\mu\|/\|\sigma_{n-p+1}-\mu\|$, where $p$ is the dimension of the desired right singular subspace of $[A,b]$ and $\mu$ is a shift for the iteration. If this ratio is not {sufficiently small}, the inverse iteration may converge very slowly. To overcome slow convergence, Chebyshev polynomials can be used, and thus, the method becomes (inverse) Chebyshev iteration. For detailed information, see \cite{v:91}. TLS problems can also have multiple right-hand sides. In this case, the block Golub-Kahan bidiagonalization procedure proposed in \cite{hps:13} can also be used to solve the problem.

  In this work, we focus on the use of Rayleigh quotient iteration (RQI) to solve the TLS problem, which is the approach advocated by Bj\"{o}rck et al. \cite{bhm:00} for large-scale problems. We introduce a mixed precision variant of the RQI-PCGTLS algorithm developed in \cite{bhm:00}. Our approach potentially decreases the computational cost of RQI-PCGTLS by using up to three different precisions in the algorithm. Using the ideas in \cite{bhm:00} and \cite{ch:22}, we discuss the convergence and accuracy of our algorithm. We derive two theoretical constraints on the precision that can be used for the construction of the preconditioner for within the inner solver. Our numerical experiments suggest that our approach can attain the same accuracy as RQI-PCGTLS, potentially with a convergence delay. Our performance modeling shows, however, that we can achieve up to $4\times$ speedup in the ideal case. Thus despite a delay in convergence, the mixed precision approach may still provide a faster time-to-solution with less computational cost.  

  We begin the next section by providing background information related to the existing RQI-based method, RQI-PCGTLS. In Section \ref{sec:rqipcgtlsmp}, we introduce our mixed precision approach, RQI-PCGTLS-MP, and present our analysis and performance modeling. In Section \ref{sec:results}, we perform numerical experiments on several small-scale matrices that are widely used in the literature for TLS problems. We finally provide our concluding remarks in Section \ref{sec:conclusion}.

  \section{Rayleigh quotient iteration with preconditioned conjugate gradient method for TLS problems (RQI-PCGTLS)}\label{sec:rqipcgtls}
  Rayleigh quotient iteration is equivalent to inverse iteration with a shift $\rho$, where \[\rho(x) = \dfrac{x^TBx}{x^Tx}\] is called the computed Rayleigh quotient of the matrix $B\in \mathbb{R}^{m\times m}$. The aim of RQI is to use the Rayleigh quotient to estimate the eigenvalue $\rho(x)$ corresponding to a given eigenvector $x$ of $B$. The computed eigenvalue is iteratively improved within the algorithm to increase the rate of convergence of inverse iteration at each step. The general scheme of RQI is given in Algorithm
\ref{rqi}.
\begin{algorithm}[htbp!]
	\caption{RQI \label{rqi}} 
	\begin{algorithmic}[1]
		\Require{Given $x_0$ with $\|x_0\|_2=1$}
		\State{$\rho_0 = (x_0)^TBx_0$} 
		\For{k = 1,2,$\ldots$}    
		\State{Solve $(B-\rho_{k-1}I)\omega=x_{k-1}$}
		\State{$x_k=\omega/\|\omega\|_2$}
		\State{$\rho_k=(x_k)^TBx_k$}
		\EndFor
	\end{algorithmic}
\end{algorithm}

  To deal with large sparse systems, Bj\"{o}rck presented an approach using RQI to solve TLS problems in \cite{b:97}. This approach solves the eigenvalue problem 
  \begin{equation}
  {\begin{bmatrix}
			A^TA &\quad A^Tb\\
			b^TA &\quad b^Tb
		\end{bmatrix}}
		\begin{bmatrix}
			x\\
			-1
		\end{bmatrix}=\lambda \begin{bmatrix}
			x\\
			-1
		\end{bmatrix}
  \label{eq:RQIsystem}
  \end{equation}
  to find $x = x_{TLS}$, where $\lambda=\sigma^2_{n+1}$. The problem can also be written as \[\begin{bmatrix}
			A^T\\
			b^T
		\end{bmatrix}(-r)=\lambda \begin{bmatrix}
			x\\
			-1
		\end{bmatrix},\] where $r = b-Ax$ is the residual. The algorithm uses Newton's method to solve 
  \begin{equation} \label{newton}
      \begin{bmatrix}
		f(x,\lambda)\\
		g(x,\lambda)
	\end{bmatrix}=\begin{bmatrix}
		-A^Tr-\lambda x\\
		-b^Tr+\lambda
	\end{bmatrix}=\begin{bmatrix}
		0\\0
	\end{bmatrix}.
  \end{equation}

 The RQI algorithm for the augmented matrix in \eqref{eq:RQIsystem} is given in Algorithm \ref{alg:rqi}. For the augmented matrix, if RQI converges to the TLS solution in $n$ steps, the Rayleigh quotient $\rho$ will be an estimate for $\sigma_{n+1}^2$. To ensure the convergence of RQI, the algorithm does one inverse iteration before starting RQI in lines \ref{line1}-\ref{line2} \cite{s:88}. Within iteration loop $k$, the Rayleigh quotient $\rho_k = \sigma_k^2=(r_k^Tr_k)/(x_k^Tx_k+1)$ is computed in line \ref{line3}. The next approximation $x_{k+1}$ is computed by solving the system \[\begin{bmatrix}
		J_k & A^Tb\\
		b^TA & \eta_k
	\end{bmatrix} \begin{bmatrix}
		x_{k+1}\\
		-1
	\end{bmatrix} = \beta_k \begin{bmatrix}
		x_k\\
		-1
	\end{bmatrix},\]
 where $J_k = A^TA-\rho_kI,$ $\eta_k = b^Tb-\rho_k,$ and $\beta_k$ is the scaling factor at the $k$-th iteration. According to \cite{bhm:00}, if $J_k$ is positive definite, the solution can be obtained by block Gaussian elimination, which leads to 
 \[
\begin{bmatrix}
		J_k & A^Tb\\
		0 & \tau_k
	\end{bmatrix} \begin{bmatrix}
		x_{k+1}\\
		-1
	\end{bmatrix} = \beta_k \begin{bmatrix}
		x_k\\
		-(z_k)^Tx_k-1
	\end{bmatrix},
 \] where $J_kz_k=A^Tb$, and $\tau_k=b^T(b-Az_k)-\rho_k$.

 We can also compute $\tau_k$ using $f_k$ and $g_k$ in each iteration using the formulas
\begin{equation}\label{compute}
    J_k\omega_k=-f_k, \quad z_k = x_k+\omega_k, \quad \tau_k = (z_k)^Tf_k-g_k.
\end{equation}
 The new way of computing $\tau_k$ can be used to express (\ref{newton}) with $f_k$ and $g_k$, i.e.,
  \begin{equation*}
      \begin{bmatrix}
		f_k\\
		g_k
	\end{bmatrix}=\begin{bmatrix}
		-A^Tr_k-\rho_k x_k\\
		-b^Tr_k+\rho_k
	\end{bmatrix},
  \end{equation*}
  where $r_k=b-Ax_k$ is the residual. 
  
 \begin{algorithm}[htbp!]
	\caption{RQI for TLS \cite{bhm:00} \label{alg:rqi}} 
	\begin{algorithmic}[1]
		%\Require{$n \times n$ matrix $A$; right-hand side $b$; maximum number of refinement steps $i_{max}$.}
		%\Ensure{Approximate solution $x_{i+1}$ to $Ax = b$.}
		\State{$x_0 = x_{LS}$} \label{line1}
		\State{$r_0  = b-Ax_0 $}
            \State{$\sigma_0 ^2 = r_0 ^Tr_0 /(x_0 ^Tx_0 +1)$}
            \State{Solve $A^TAu_0  = x_0 $}
            \State{$x_1=x_0 +\sigma_0 ^2 u_0 $}\label{line2}
		\For{{$k$}=1,2,$\ldots$}    
		\State{$r_k = b-Ax_k$}
		\State{$\sigma_k^2 = r_k^Tr_k/(1+x_k^Tx_k)$} \label{line3}
		\State{$f_k=-A^Tr_k-\sigma_k^2x$}
            \State{$g_k = -b^Tr_k+\sigma_k^2$}
            \State{Solve $(A^TA-\sigma_k^2I)\omega_k = -f_k$}\label{solve1}
            \State{$z_k= x_k+\omega_k$}
            \State{$\beta_k = (z_k^Tf_k-g_k)/(z_k^Tx_k+1)$}
            \State{Solve $(A^TA-\sigma_k^2I)u_k = x_k$}\label{solve2}
            \State{$x_{k+1}=z_k+\beta_k u_k$}
		\EndFor
	\end{algorithmic}
\end{algorithm}

Performing one Rayleigh quotient iteration requires solving two linear systems (lines \ref{solve1} and \ref{solve2}). Since the coefficient matrix $A^TA-\sigma^2I$ is symmetric and positive definite, the conjugate gradient method can be applied to the problem. In \cite{bhm:00}, the RQI algorithm is used with a preconditioned conjugate gradient algorithm called PCGTLS to approximately  solve these two systems using a fixed number of steps with the Cholesky factor of $A^TA$ used as the preconditioner. Bj\"{o}rck et al. call this approach RQI-PCGTLS \cite{bhm:00}. Note that we need to form $A^TA$ and compute its Cholesky factor $R$ only once. The PCGTLS algorithm is given in Algorithm \ref{alg:pcgtls}. The algorithm performs $k+1$ PCG steps in the $k$-th Rayleigh quotient iteration; the number of PCG steps is discussed further in \cite{bhm:00}.
 \begin{algorithm}[htbp!]
	\caption{{$l$-step} PCGTLS for $(A^TA-\sigma^2I)\omega = f$ \cite{bhm:00} \label{alg:pcgtls}} 
	\begin{algorithmic}[1]
		%\Require{$n \times n$ matrix $A$; right-hand side $b$; maximum number of refinement steps $i_{max}$.}
		%\Ensure{Approximate solution $x_{i+1}$ to $Ax = b$.}
		\State{Initialize $\omega_0=0$, $p_0=s_0=R^{-T}f$, $\eta_0 = \|s_0\|^2_2$}\label{p1}
		\For{$j=0,1,\ldots,l$, while $\delta_j\neq 0$}    
		\State{$q_j=R^{-1}p_j$}\label{p2}
		\State{$\delta_j = \|p_j\|^2_2-\sigma^2\|q_j\|^2_2$}\label{pc1}
		\State{$\alpha_j = \eta_j/\delta_j$}\label{pc2}
            \State{$\omega_{j+1}=\omega_j+\alpha_jq_j$}\label{pc3}
            \State{$q_j=R^{-T}q_j$}\label{p3}
            \State{$s_{j+1}=s_j-\alpha_j(p_j-\sigma^2q_j)$}\label{pc4}
            \State{$\eta_{j+1}=\|s_{j+1}\|^2_2$}\label{pc5}
            \State{$\beta_j=\eta_{j+1}/\eta_j$}\label{pc6}
            \State{$p_{j+1}=s_{j+1}+\beta_jp_j$}\label{pc7}
		\EndFor
	\end{algorithmic}
\end{algorithm}

In Algorithm \ref{alg:pcgtls}, the Cholesky factor $R$ of $A^TA$ is assumed to be computed exactly and used to explicitly form the preconditioned matrix $\Tilde{C}= I-\sigma^2R^{-T}R^{-1}$. Thus, the matrix $A$ is never used in PCGTLS except in forming the right-hand side vector $f$, which reduces the computational cost of the algorithm significantly. Therefore, forming the matrix $A^TA$ and computing its Cholesky factor can be considered the most expensive parts of RQI-PCGTLS. This is especially true if $A$ is not skinny, i.e., if $n$ is close to $m$, since forming $A^TA$ has cost $\mathcal{O}(n^2m)$, and Cholesky factorization has cost $\mathcal{O}(n^3)$.

 \section{Mixed precision RQI-PCGTLS (RQI-PCGTLS-MP)}  \label{sec:rqipcgtlsmp}
 One can use several techniques to reduce the computational cost of an algorithm. A popular approach, motivated by the recent emergence of commercially-available low-precision hardware, is the development of mixed precision algorithms, in which lower precision is used for the computationally dominant parts. However, one must use lower precision with caution, since this can cause loss of accuracy, convergence delay, and make the algorithm more susceptible to overflow and underflow. 
 
Our goal is thus to use different precisions in different parts of the algorithm to safely balance accuracy, stability, and cost.  In \cite{fs:07}, the authors state that the inverse iteration need not be performed exactly to be convergent. We therefore expect that mixed precision can be used within the RQI-PCGTLS algorithm. Moreover, assume that RQI solves the eigenvalue problem $(A-\Theta I)\omega=z$ for the eigenpair $(\Theta,\omega)$ where $z$ is an approximate eigenvector. The authors in \cite{se:02} proved that this is mathematically equivalent to Newton's method on the unit sphere, i.e., solving $\Pi (A-\Theta I)\Pi d = -r$, where $d$ is a correction with $z^*d=0$, $\Pi=I-zz^*$, and $r=\Pi Az=Az-\Theta z$, if a Galerkin method is used to solve both systems. Thus, motivated by the mixed precision inexact Newton method introduced in \cite{k:22}, we developed a mixed precision variant of RQI-PCGTLS, which we call RQI-PCGTLS-MP. We note that we expect that the analyses of RQI in \cite{fs:07,se:02} and Newton in \cite{hm:22,k:22,t:01} can be extended to give a complete analysis of RQI-PCGTLS-MP; we leave this as future work. 
 
 Our mixed precision approach, RQI-PCGTLS-MP, is presented in Algorithm \ref{alg:rqi_mp}. RQI-PCGTLS-MP uses three different precisions: $u$ is the working precision to store data and solutions in the RQI-PCGTLS-MP algorithm, $u_p$ is the working precision in the PCGTLS algorithm, and $u_q$ is the precision for computing the QR factors of $A$. Throughout this work, we assume that $u\geq u_p\geq u_q$.
 
 The authors of \cite{bhm:00} noted that one could use the R-factor from the QR decomposition of $A$ instead of the Cholesky factor of $A^TA$. Since $A$ is rectangular, the QR decomposition of the matrix can be written \[A = [Q_1 \quad Q_2]\begin{bmatrix}
     R\\
     0
 \end{bmatrix}.\] From \cite{b:96}, we know that if $A$ is of full rank, the R-factor in the QR decomposition of $A$ is equal to the Cholesky factor of $A^TA$. However, the use of the Cholesky factor of $A^TA$ results in a squaring of $\kappa(A)$, which is undesirable from a numerical perspective, especially if we hope to use low precision for this computation. Thus to enable the use of lower precision for more ill-conditioned systems, we use the R-factor from the Householder QR factorization of $A$ instead of the Cholesky factorization of $A^TA$ within RQI-PCGTLS-MP. We note that we have, however, observed experimentally that in some cases using the Cholesky factor of $A^TA$ over the R-factor of $A$ may be better in terms of RQI convergence rate. Determining the mechanism for these observations is left as future work.

\begin{algorithm}[htbp!]
	\caption{RQI-PCGTLS-MP \label{alg:rqi_mp}} 
	\begin{algorithmic}[1]
		%\Require{$n \times n$ matrix $A$; right-hand side $b$; maximum number of refinement steps $i_{max}$.}
		%\Ensure{Approximate solution $x_{i+1}$ to $Ax = b$.}
		\State{$x_0 = x_{LS}$}
		\State{$r_0 = b-Ax_0$}{\hspace*{\fill} {in precision $u$}}\label{mp1}
            \Statex{{(\textit{$2mn+m$ ops})}}
            \State{$\sigma_0^2 = r_0^Tr_0/(1+x_0^Tx_0)$}{\hspace*{\fill} {in precision $u$}}\label{mp2}
            \Statex{{(\textit{$2m+2n-2$ ops})}}
            \State{Compute $A=Q[R\quad 0]^T$}{\hspace*{\fill} {in precision $u_q$}}\label{mp3}
            \Statex{{(\textit{$2mn^2-2n^3/3$ ops})}}
            \State{Solve $A^TAu_0 = x_0$}{\hspace*{\fill} {in precision $u$}}\label{mp4}
            \Statex{{(\textit{$2n^2$ ops})}}
            \State{$x_1=x_0+\sigma_0^2 u_0$}{\hspace*{\fill} {in precision $u$}}\label{mp5}
            \Statex{{(\textit{$2n$ ops})}}
		\For{k = 1,2,$\ldots$}    
		\State{$r_k = b-Ax_k$}{\hspace*{\fill} {in precision $u$}}\label{mp6}
            \Statex{\quad{(\textit{$2mn+m$ ops})}}
		\State{$\sigma_k^2 = r_k^Tr_k/(1+x_k^Tx_k)$}{\hspace*{\fill} {in precision $u$}}\label{mp7}
            \Statex{\quad{(\textit{$2m+2n-2$ ops})}}
		\State{$f_k=-A^Tr_k-\sigma_k^2x_k$}{\hspace*{\fill} {in precision $u$}}\label{mp8}
            \Statex{\quad{(\textit{$2mn+2n$ ops})}}
            \State{$g_k = -b^Tr_k+\sigma_k^2$}{\hspace*{\fill} {in precision $u$}}\label{mp9}
            \Statex{\quad{(\textit{$2m-1$ ops})}}
            \State{Solve $(A^TA-\sigma_k^2I)\omega_k = -f_k$ via PCGTLS ($\leq k+1$ iterations)}{\hspace*{\fill} {in precision $u_p$}}\label{algline1}
        \Statex{\quad{(\textit{$2n^2+14n-3$ ops per PCGTLS iteration $+n^2+2n-1$ ops})}}
            \State{$z_k= x_k+\omega_k$}{\hspace*{\fill} {in precision $u$}}\label{mp10}
        \Statex{\quad{(\textit{$n$ ops})}}
            \State{$\beta_k = (z_k^Tf_k-g_k)/(z_k^Tx_k+1)$}{\hspace*{\fill} {in precision $u$}}\label{mp11}
        \Statex{\quad{(\textit{$4n-2$ ops})}}
            \State{Solve $(A^TA-\sigma_k^2I)u_k = x_k$ via PCGTLS ($\leq k+1$ iterations)}{\hspace*{\fill} {in precision $u_p$}}\label{algline2}
    \Statex{\quad{(\textit{$2n^2+14n-3$ ops per PCGTLS iteration $+n^2+2n-1$ ops})}}
            \State{$x_{k+1}=z_k+\beta_k u_k$}{\hspace*{\fill} {in precision $u$}}\label{mp12}
        \Statex{\quad{(\textit{$2n$ ops})}}
		\EndFor
	\end{algorithmic}
\end{algorithm}

To terminate the Rayleigh quotient iterations in Algorithm \ref{alg:rqi_mp}, we use the termination criteria given in \cite{bhm:00}. The authors in \cite{bhm:00} discuss that the convergence is directly related to the behavior of the normalized residual norm, 
\begin{equation}
    \psi_k = \left(\frac{\|f_k\|^2_2+g_k^2}{\|x_k\|_2^2+1}\right)^{1/2}.
    \label{eq:psik}
\end{equation}
The method is converging in exact arithmetic as long as $\psi_{k+1}\leq\psi_k$. Hence, when $\psi_{k+1}>\psi_k$, we terminate the iterations. 

Assuming exact arithmetic, there are several necessary conditions for the RQI-PCGTLS-MP algorithm (and thus also the RQI-PCGTLS algorithm) to converge. 
The first condition guarantees that the TLS problem has a unique solution, $x_{TLS}$, under exact arithmetic. If $\sigma'_n>\sigma_{n+1}$, then $x_{TLS}$ is unique \cite{b:96}. The second condition is based on the assumptions required for the PCGTLS algorithm. If $A^TA-\sigma_k^2I$ is symmetric and positive definite for all $\sigma_k$, then PCGTLS converges to an approximate solution in exact arithmetic \cite{bhm:00}. 

\subsection{Constraints on factorization precision}

In the case of very large $A$, the QR factorization (or Cholesky factorization) to construct the preconditioner is expected to be the most expensive part of the computation. Thus, we want to use as low a precision as possible for $u_q$. There are two potential sources from which a constraint on $u_q$ arises. 

First, we want {the} computed $\hat{R}$ factor to be nonsingular, since we will use it as a preconditioner. Note that here and in the remainder of the paper we use hats to denote quantities computed in finite precision. 
The precision required to guarantee the nonsingularity of $\hat{R}$ can be found in, e.g., \cite[Theorem 19.4]{h:02}. Suppose the Householder QR factorization of $A$ is computed in precision $u_q$ so that $A+\Delta A=\tilde{Q}_1\hat{R}$, where $\hat{R}\in \mathbb{R}^{n\times n}$ and $\tilde{Q}_1\in\mathbb{R}^{m\times n}$ has orthonormal columns. If 
\[
u_q<\dfrac{1}{cmn^{3/2}\kappa(A)},
\] 
for a small integer constant $c$, then $\hat{R}$ is guaranteed to be non-singular. We note that a probabilistic approach to rounding error analysis can loosen this bound by reducing the dimensional factors by their square roots. {In the probabilistic model, the rounding errors $\gamma_1, \gamma_2, \ldots$ satisfying $fl(x\enspace \text{op}\enspace y) = (x\enspace \text{op}\enspace y)(1+\delta)$ with $|\delta|\leq u$ and op $=+,-,*,/$ are mean independent random variables of mean zero.} The analysis in \cite[Theorem 4.4]{ch:22} { proves that under the probabilistic model of rounding errors, applying the matrix concentration equality given in }\cite[Theorem 3.2]{ch:22}{ results in a bound on the error of Householder QR in which dimensional factors may be replaced by their square roots, thus giving the constraint}

\begin{equation*}
    u_q<\dfrac{1}{cm^{1/2}n^{3/4}\kappa(A)}.
\end{equation*}

As a rough heuristic, this means that we must have 
\begin{equation}\label{constr1}
    u_q \lesssim \kappa(A)^{-1}.
\end{equation}

A similar constraint can be obtained for the case where the Cholesky factorization of $fl(A^TA)$ is used to construct the preconditioner. From Wilkinson's result on Cholesky factorization in \cite{w:71}, we say the Cholesky factorization of $fl(A^TA)$ succeeds in precision $u_q$ if 
\[
u_q\leq \dfrac{1}{20n^{3/2}\kappa(A)^2},
\]
or as a rough heuristic, 
\[
u_q \lesssim \kappa(A)^{-2}.
\]
Further, using the results in \cite{d:89} and \cite{hp:21},
if we introduce the scaled matrix $H$ given by 
\[
H = D^{-1}fl(A^TA)D^{-1}, 
\]
where $D$ is a diagonal matrix with the square root of the diagonal elements of $A^TA$ on the diagonal, we can say that if 
\[u_q<\dfrac{\lambda_{min}(H)}{(2\lambda_{min}(H)+n)(n+1)},\] 
where $\lambda_{min}(H)$ is the minimum eigenvalue of $H$, then the Cholesky factorization of $fl(A^TA)$ with two-sided diagonal scaling, run on $fl(A^TA)$ computed in precision $u_q$, will succeed.

Note that this result assumes that we can apply this scaling without incurring additional rounding error. For a dynamic approach to scaling and shifting for low-precision Cholesky factorization along with a complete analysis, see \cite{hp:21}.  In our numerical experiment involving Cholesky factorization in Section \ref{sec:results}, we will apply a two-sided scaling for equilibration before Cholesky factorization to avoid overflow and underflow, which may result from rounding higher precision numbers to a lower precision.

Our second constraint on $u_q$ will come from the requirement that the preconditioned system constructed using the computed $\hat{R}$ must remain positive definite.  
Assume we are computing $R$ via Householder QR factorization. Then we have 
\begin{equation}
A+ \Delta A = Q_1\hat{R} \quad \text{with} \quad \|\Delta A\|_F \leq \Tilde{\gamma}_{mn}^q \|A\|_F, \quad \text{where} \quad \Tilde{\gamma}_{mn}^q=\dfrac{cmnu_q}{1-cmnu_q}
\label{eq:deltA}
\end{equation}
for some constant $c$, where the subscript $F$ denotes the Frobenius norm. Note that if $A=Q_1R$, we have 
\begin{equation}
\| \tilde{Q}_1-Q_1 \|_F \leq \sqrt{m}\tilde{\gamma}_{mn}^q;
\label{eq:Qdiff}
\end{equation}
see, e.g., \cite[Equation (2.3)]{chp:20}.

When we take the inexactness of $\hat{R}$ due to computation in precision $u_q$ into account in determining the preconditioned coefficient matrix $\hat{C}$, we have 
\begin{align}
\hat{C} & =\hat{R}^{-T}(A^TA-\sigma^2I)\hat{R}^{-1}\nonumber\\
 & = \hat{R}^{-T}A^TA\hat{R}^{-1}-\sigma^2\hat{R}^{-T}\hat{R}^{-1} \nonumber\\
 & = (\tilde{Q}_1^T-\hat{R}^{-T}\Delta A^T)(\tilde{Q}_1-\Delta A\hat{R}^{-1})-\sigma^2\hat{R}^{-T}\hat{R}^{-1} \nonumber\\
 & \approx I - \tilde{Q}_1^T\Delta A \hat{R}^{-1}-\hat{R}^{-T}\Delta A^T\tilde{Q}_1 - \sigma^2\hat{R}^{-T}\hat{R}^{-1}.\label{eq:Chat}
\end{align}
We must now investigate the last term on the right-hand side above. We can write (using a first order approximation)
\begin{align}
    \sigma^2\hat{R}^{-T}\hat{R}^{-1} &= \sigma^2 \left((\hat{R}^T - R^T) + R^T \right)^{-1}\left((\hat{R} - R) + R \right)^{-1}  \nonumber\\
    &= \sigma^2 \left( R^T \left(I-R^{-T}(R^T - \hat{R}^T)\right) \right)^{-1}\left( R \left(I-R^{-1}(R - \hat{R})\right) \right)^{-1} \nonumber\\
    &\approx \sigma^2 \left(R^{-T} + R^{-T}(R^T-\hat{R}^T)R^{-T} \right)\left(R^{-1} + R^{-1}(R-\hat{R})R^{-1} \right) \nonumber \\
    &\approx \sigma^2 \left(R^{-T}R^{-1} + R^{-T}R^{-1}(R-\hat{R})R^{-1} + R^{-T}(R^T-\hat{R}^T)R^{-T}R^{-1} \right). \nonumber
\end{align}
Using that 
\[
\hat{R}-R = (\tilde{Q}_1 - Q_1)^T A + \tilde{Q}_1^T\Delta A,
\]
the above becomes
\begin{align*}
 \sigma^2\hat{R}^{-T}\hat{R}^{-1} &\approx \sigma^2 \Big(R^{-T}R^{-1} + R^{-T}R^{-1}\big( (\tilde{Q}_1 - Q_1)^T A + \tilde{Q}_1^T\Delta A \big) R^{-1} \\
 &\phantom{\approx} + R^{-T}\big( A^T(\tilde{Q}_1 - Q_1) + \Delta A^T \tilde{Q}_1\big)R^{-T}R^{-1} \Big). 
\end{align*}
Substituting this into {\eqref{eq:Chat}}, we have
\begin{align*}
    \hat{C} &\approx I - \sigma^2 R^{-T}R^{-1} - \tilde{Q}_1^T\Delta A \hat{R}^{-1}-\hat{R}^{-T}\Delta A^T\tilde{Q}_1 + \sigma^2 R^{-T}R^{-1}\big( (\tilde{Q}_1 - Q_1)^T A + \tilde{Q}_1^T\Delta A \big) R^{-1}  \\
    &\phantom{\approx}+ \sigma^2 R^{-T}\big( A^T(\tilde{Q}_1 - Q_1) + \Delta A^T \tilde{Q}_1\big)R^{-T}R^{-1} \\
    &\equiv\tilde{C}+\underline\Delta,
\end{align*}
where $\tilde{C}$ is the preconditioned system with exact $R$, i.e., $\tilde{C}=I-\sigma^2R^{-T}R^{-1}$. 

We now want to bound the norm of $\underline\Delta$. Using {the definition of $\underline\Delta$ above, along with} {\eqref{eq:deltA}}{ and }{\eqref{eq:Qdiff}}{, and ignoring terms of order $O(u_q^2)$, we have}
\begin{align*}
    \|\underline\Delta\|_2 &\leq 2\|\Delta A\|_2 \|\hat{R}^{-1}\|_2 + 2 \sigma^2 \|\hat{R}^{-1}\|_2^{3} \left( \|\tilde{Q}_1 - Q_1\|_2 \|A\|_2 + \|\Delta A\|_2  \right)\\
    &\leq \tilde{\gamma}_{mn}^q \kappa_F(A) + \sigma^2 \|A^{-1}\|_2^3 \cdot \sqrt{m}\tilde{\gamma}_{mn}^q \|A\|_2 \\
    &\leq \left(1 + \sigma^2 \|A^{-1}\|_2^2 \right) \sqrt{m} \tilde{\gamma}_{mn}^q \kappa_F(A),
\end{align*}
where $\kappa_F(A)$ denotes the Frobenius norm condition number. 
Assuming that $\sigma =\sigma_{n+1}$, we have $\sigma \leq 1/\|A^{-1}\|_2$, and thus 
\[
\|\underline\Delta\|_2 \leq \sqrt{m} \tilde{\gamma}_{mn}^q \kappa_F(A). 
\]

Using eigenvalue perturbation theory, we can say how far the eigenvalues of $\hat{C}$ are from $\Tilde{C}$ for the case $\sigma = \sigma_{n+1}$. From (4.4) in \cite{bhm:00}, we know that in this case the eigenvalues of $\tilde{C}$ lie in the interval \[\left[1-\dfrac{\sigma_{n+1}^2}{(\sigma'_n)^2},1\right],\]
and so we expect the exact preconditioner to work very well unless $\sigma_{n+1}^2 \approx (\sigma'_n)^2$, {in which case the condition number may be large and thus fast convergence of CG is not guaranteed.}

For $\hat{C}$ with {\emph{inexact}} $\hat{R}$, we know the eigenvalues must lie in the interval \[\left[1-\dfrac{\sigma_{n+1}^2}{(\sigma'_n)^2}-\|\underline\Delta\|_2,1+\|\underline\Delta\|_2\right].\]
This tells us two things. First, for the computed $\hat{R}$, we expect the preconditioner to work well unless \[\sigma_{n+1}^2 \approx (1+\|\underline\Delta\|_2)(\sigma'_n)^2.\]
Second, we must have that 
\[1-\dfrac{\sigma_{n+1}^2}{(\sigma'_n)^2}-\|\underline\Delta\|_2>0,
\] 
which implies
\[
\|\underline\Delta\|_2<1-\dfrac{\sigma_{n+1}^2}{(\sigma'_n)^2},
\]
or else we cannot guarantee the positive definiteness of the preconditioned system. 

As a rough heuristic, this means that we want 
\[
u_q \kappa_F(A)\lesssim 1-\dfrac{\sigma_{n+1}^2}{(\sigma'_n)^2}.
\]
In other words, we must have 
\begin{equation}\label{constr2}
    u_q\lesssim {\kappa_F(A)^{-1}}\left(1-\dfrac{\sigma_{n+1}^2}{(\sigma'_n)^2}\right).
\end{equation}

Note that this means that for TLS problems, in contrast to standard least squares problems, our choice of precision should depend not only on $A$ but also on the right-hand side $b$.

To compare the bounds given in {\eqref{constr1}} and {\eqref{constr2}}, we need to look at $\sigma_{n+1}^2 / (\sigma'_n)^2$. When $\sigma_{n+1}$ is close to zero, i.e., $b$ is almost in the column space of $A$, $\sigma_{n+1}^2 /(\sigma'_n)^2$ is very small, and thus the two constraints become very close to each other. Otherwise, if $\sigma_{n+1}^2 /(\sigma'_n)^2$ is close to 1, then the constraint on $u_q$ in \eqref{constr2} is much tighter. In this case, we need to use a higher precision for $u_q$.

\subsection{Performance modeling}
To evaluate to what extent the computational cost can be reduced by using the mixed precision variant with Householder QR factorization, we construct a performance model which takes as input the size $(m,n)$ of $A$ (with $m\geq n$) and the number $r$ of the Rayleigh quotient iterations performed by RQI-PCGTLS and RQI-PCGTLS-MP. The model assumes that PCGTLS performs exactly $k+1$ iterations in the $k$-th Rayleigh quotient iteration. We also note that our performance model is based on computational complexity and thus essentially assumes an ideal model in which data movement is not dominant; this may not be realistic in practice.  

{Using the operation counts in Table }{\ref{tab:cost},} we calculate the cost of each line in Algorithm \ref{alg:rqi_mp}. Then we multiply them by a constant depending on the chosen precision. For unit roundoff of different precisions, see Table \ref{tab:eps}. Since the default working precision is fp64, we use 1 for double, 0.5 for single, and 0.25 as the constant for half precision. For instance, line \ref{mp1} in Algorithm \ref{alg:rqi_mp}, $r=b-Ax$, performs one matrix-vector product and a vector update. Therefore this line requires $2mn-m+2m=2mn+m$ operations. Since we use double precision for $u$ in our experiments, the cost of this line becomes $1\times (2mn+m)=2mn+m$. In the end, the total operation count for the RQI-PCGTLS-MP algorithm depends on the size of $A$, the precision setting, and the number $r$ of Rayleigh quotient iterations performed and can be written as

\begin{align*}
    cost(m,n,r,c,c_p,c_q) &= c (2mn+3m+4n+2n^2-2+r(4mn+{5}m+11n-5))\\
    &+ c_p (2r(n^2+2n-1)+(r^2+3r)(2n^2+14n-3))\\
    &+ c_q \left(2mn^2-\dfrac{2n^3}{3}\right),
\end{align*}
where $(c,c_p,c_q)$ are the constants used for $(u,u_p,u_q)$, respectively. For the RQI-PCGTLS algorithm, we use the same $cost$ function with $(c,c_p,c_q) = (1,1,1)$. Then we calculate the speedup of RQI-PCGTLS-MP as 
\[
speedup = \dfrac{cost \text{ of RQI-PCGTLS-MP}}{cost\text{ of RQI-PCGTLS}}.
\]

\begin{table}[]
\centering
\caption{{Cost and places of each operation used in RQI-PCGTLS-MP. In the table, we use the abbreviations QR decomp (QR decomposition), Subs (Substitutions), MatVec (Matrix-Vector multiplication $Ab$), VecMul (Inner product $b^Tc$), VecUpdate (Vector update $\lambda b+c$). }}
\resizebox{\linewidth}{!}{\begin{minipage}{\textwidth}
\label{tab:cost}
\begin{tabular}{cccc}
operation & cost & lines in Algorithm \ref{alg:rqi_mp}& lines in Algorithm \ref{alg:pcgtls}\\ \hline
QR decomp of $A^{m\times n}$ &  $2mn^2-2n^3/3$   &  \ref{mp3}  &   - \\
Forward/Backward Subs with $R^{n\times n}$ & $n^2$&  \ref{mp4} &  \ref{p1}, \ref{p2}, \ref{p3} \\
MatVec  with $A^{m\times n}, b^{n\times 1}$  &  $2mn-m$  & \ref{mp1}, \ref{mp6}, \ref{mp8}   &   -   \\
VecMul with $b^{n\times 1},c^{n\times 1}$   & $2n-1$     & \ref{mp2}, \ref{mp7}, \ref{mp9}, \ref{mp11}   & \ref{p1}, \ref{pc1}, \ref{pc5}   \\
VecUpdate with $b^{n\times 1},c^{n\times 1},\lambda\in \mathbb{R}$ & $2n$  &  \ref{mp1}, \ref{mp5}, \ref{mp6}, \ref{mp8}, \ref{mp10}, \ref{mp12}  &  \ref{pc3}, \ref{pc4}, \ref{pc7}  
\end{tabular}
\end{minipage}}
\end{table}

\begin{table}[]
	\centering
	\caption{Various IEEE precisions and their units roundoff. }
	\label{tab:eps}
	\begin{tabular}{|l|l|}
		\hline
		\multicolumn{1}{|c|}{Precision} & \multicolumn{1}{c|}{Unit Roundoff} \\ \hline
		fp16 (half) & $4.88\cdot 10^{-4}$ \\ \hline
		fp32 (single) & $5.96\cdot 10^{-8}$ \\ \hline
		fp64 (double) & $1.11\cdot 10^{-16}$ \\ \hline
		fp128 (quad) & $9.63\cdot 10^{-35}$ \\ \hline
	\end{tabular}
\end{table}

\begin{figure}
    \centering
    \includegraphics[width=0.48\linewidth]{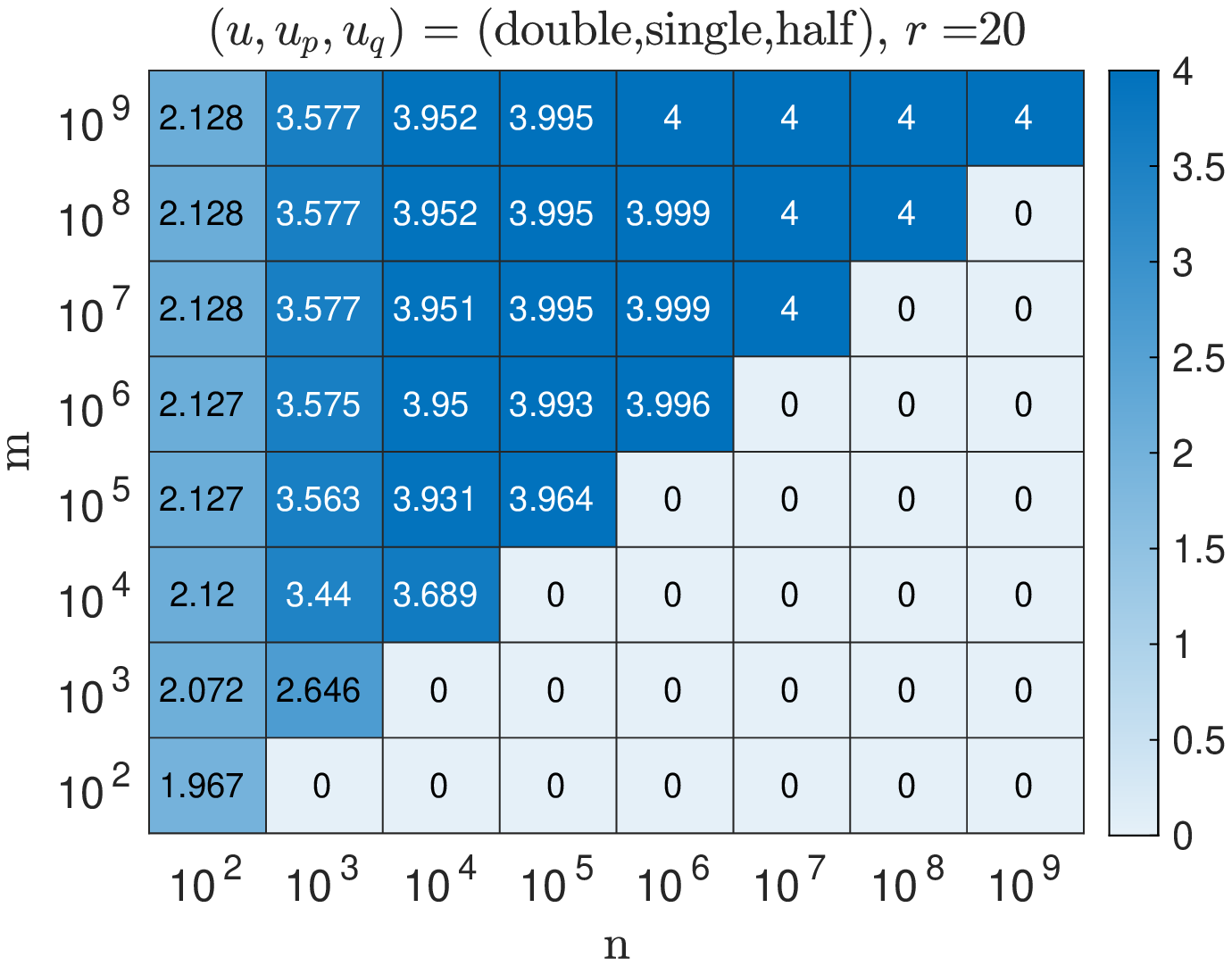}
    \vspace{0.5cm}
    \hspace{0.2cm}
    \includegraphics[width=0.48\linewidth]{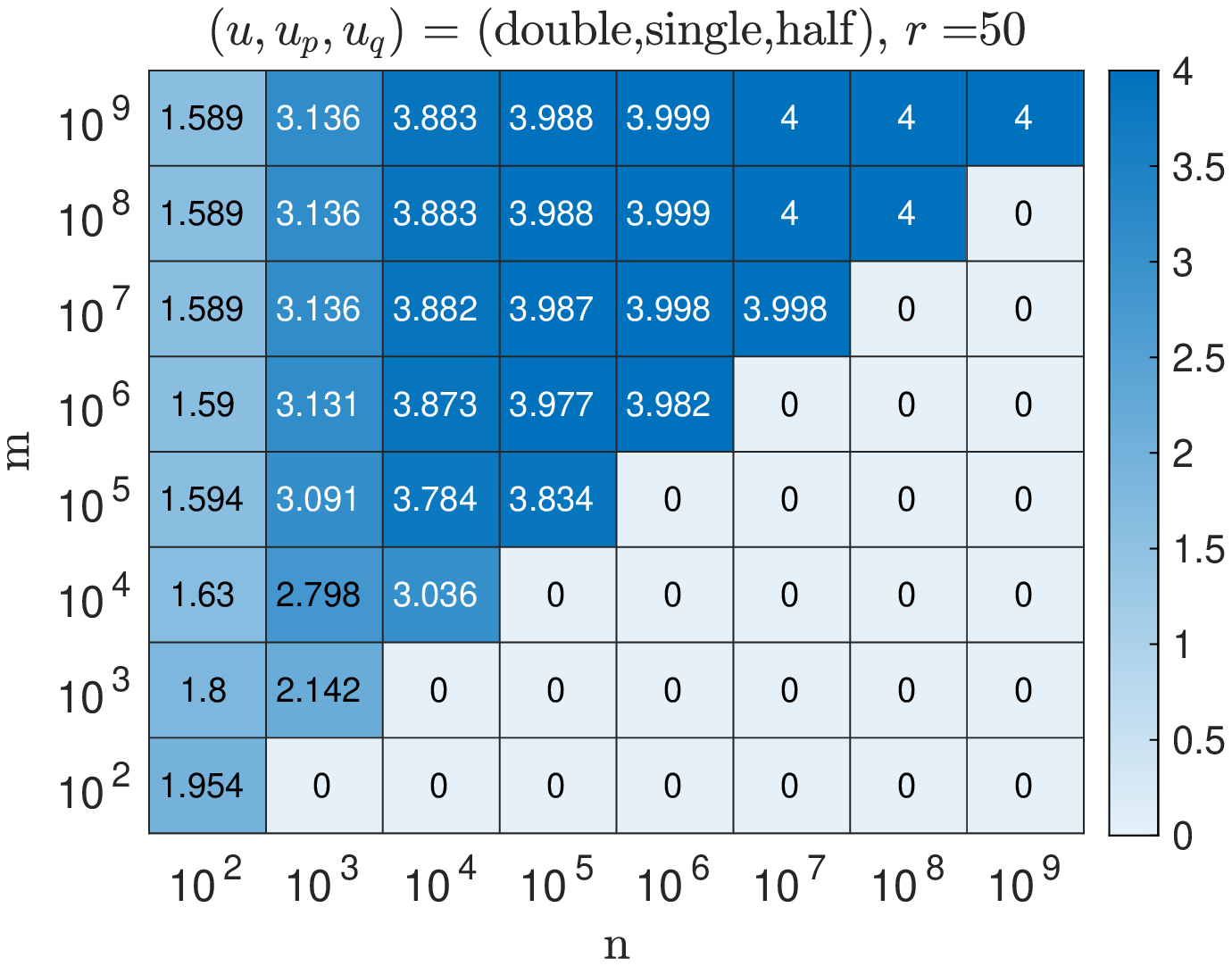}
    \vspace{0.5cm}
    \includegraphics[width=0.48\linewidth]{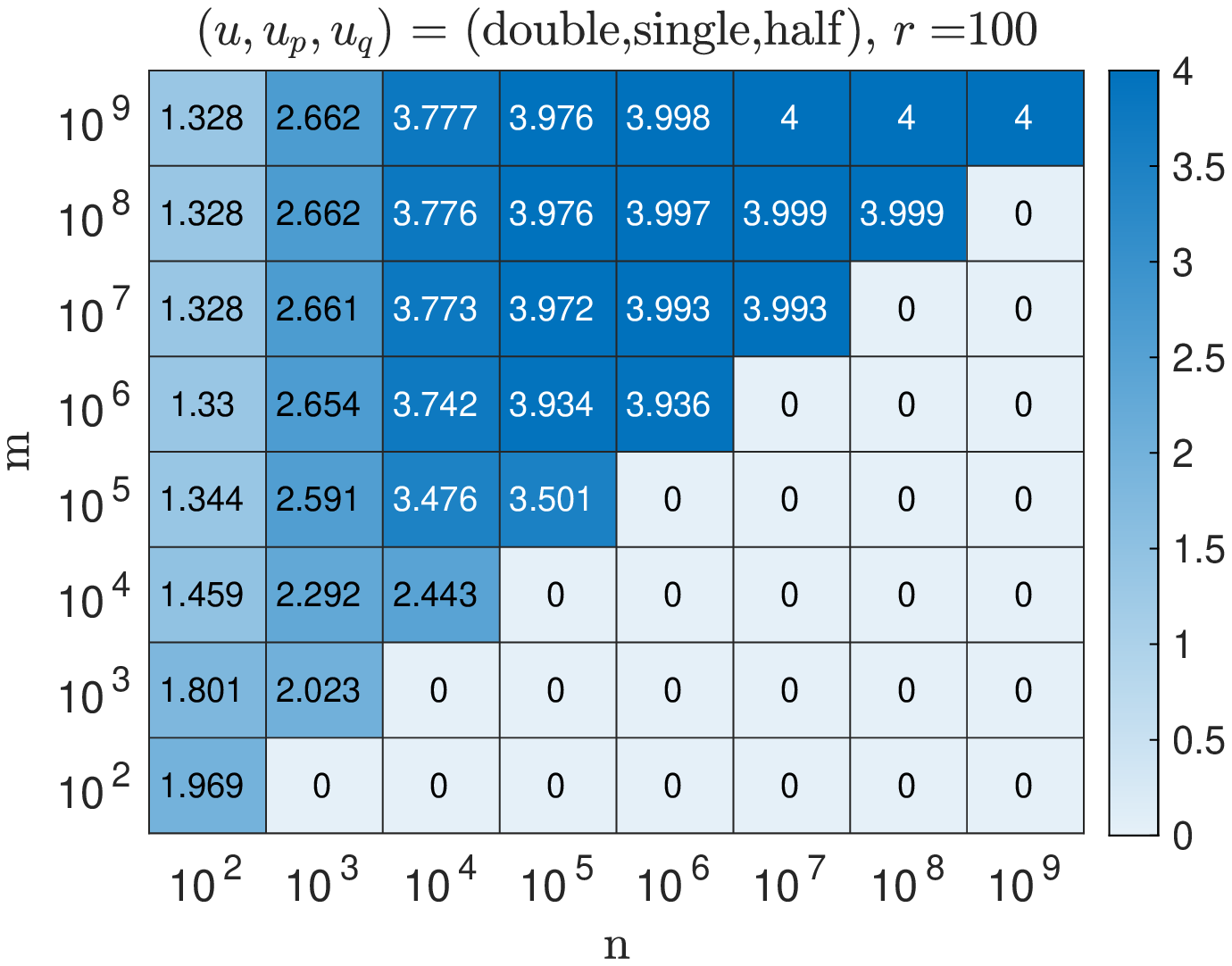}
    \hspace{0.2cm}
    \includegraphics[width=0.48\linewidth]{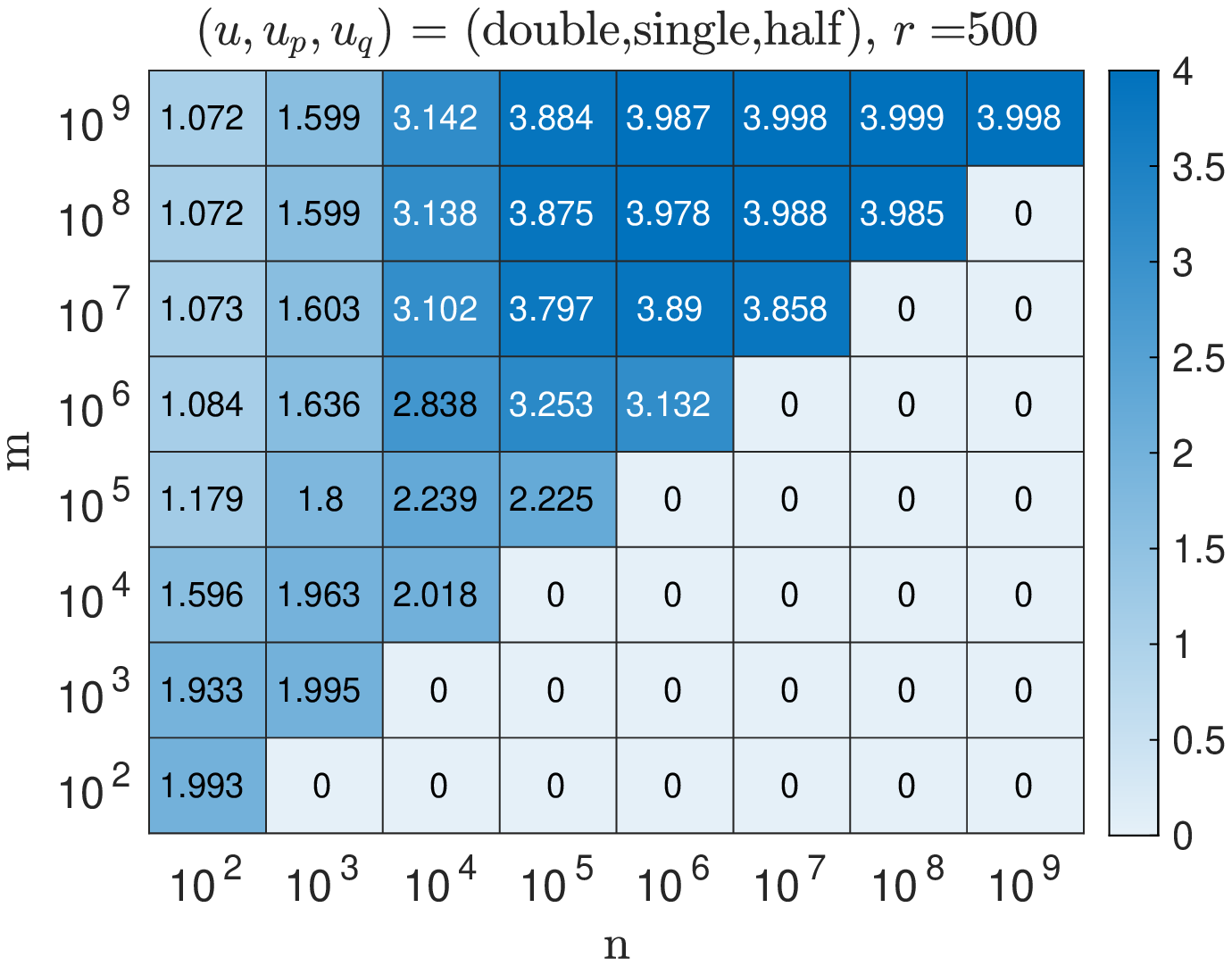}
    \includegraphics[width=0.48\linewidth]{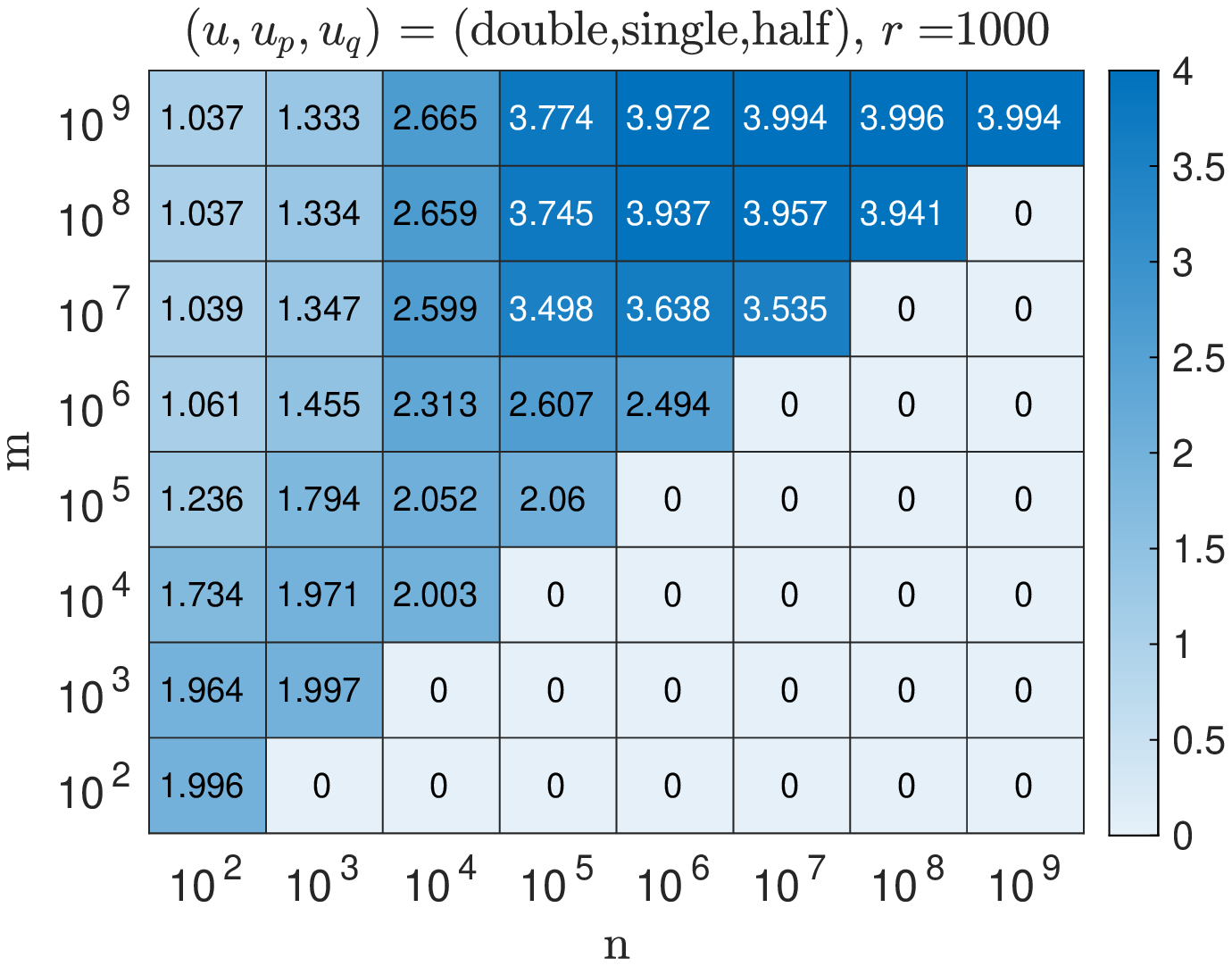}
    \hspace{0.2cm}
    \includegraphics[width=0.48\linewidth]{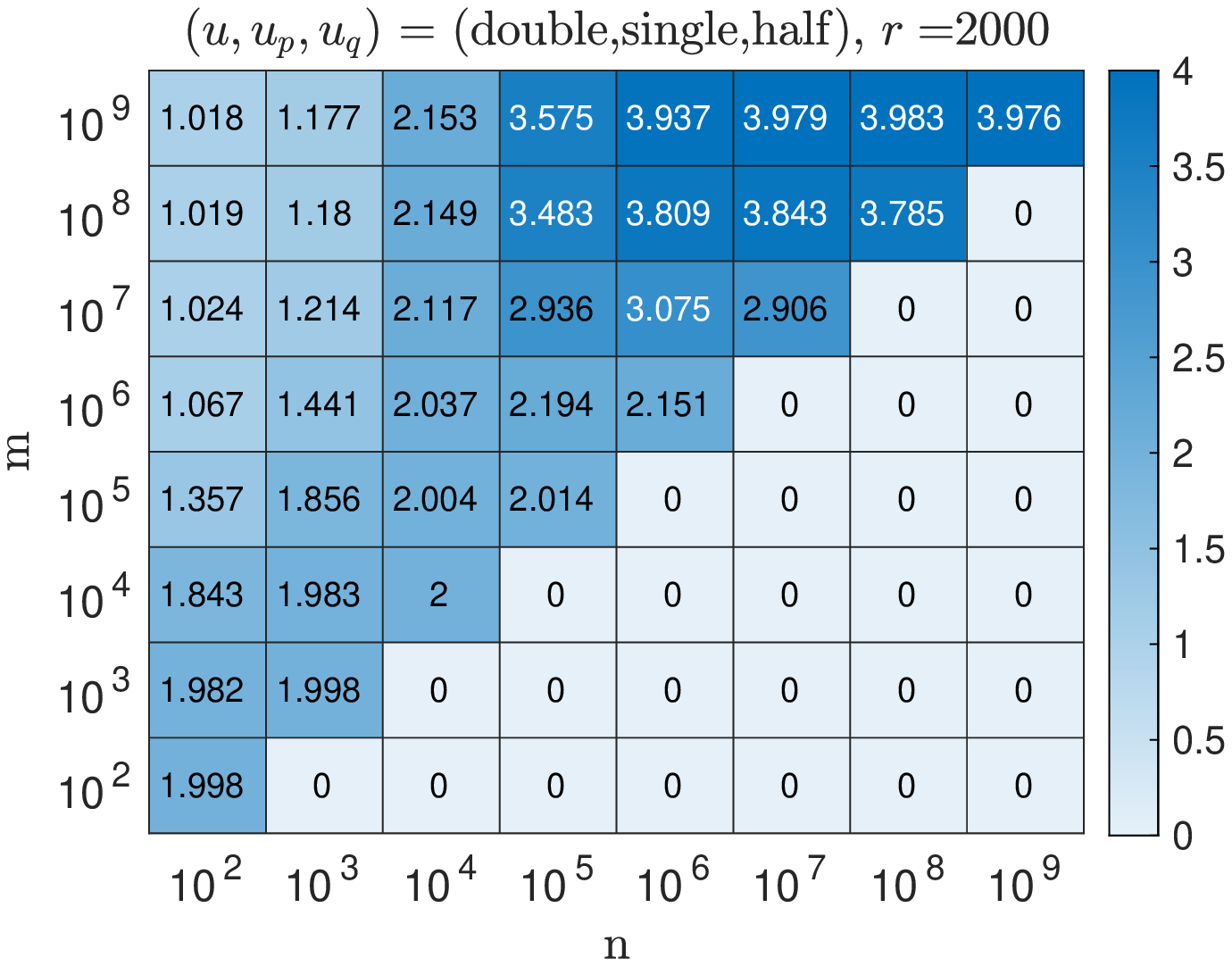}
    \caption{{Modeled speedup for $A\in \mathbb{R}^{m\times n}$ (for $m\geq n$) with $(u,u_p,u_q)$ = (double, single, half) and various numbers $r$ of Rayleigh quotient iterations performed by RQI-PCGTLS and RQI-PCGTLS-MP.}}
    \label{fig:performance}
\end{figure}

Figure \ref{fig:performance} illustrates the speedups one can obtain using the $(u,u_p,u_q)$ = (double, single, half) setting with several values of $r$. Using this mixed precision setting can provide up to $\times 4$ speedup in the best-case scenario, which comes from the use of half precision (fp16) for the dominant part of the computation. From the figure, we observe that even for a tall-skinny matrix, our mixed precision approach can perform more than $2$ times faster than the RQI-PCGTLS algorithm. We also see that we get a greater speedup the larger $n$ grows. The best-case scenario is obtained when both $m$ and $n$ are large. In these cases, we can obtain $4 \times$ modeled speedup even with $2000$ Rayleigh quotient iterations. One should note that $r$ depends on the matrix. Thus for a large-scale matrix, $r$ may be larger than $2000$. However, the speedup will never drop below 1 as long as $u$ is not higher than the precision RQI-PCGTLS uses.

	\section{Numerical experiments}\label{sec:results}
In this section, we present the numerical results comparing RQI-PCGTLS and RQI-PCGTLS-MP for solving TLS problems that are commonly used in the literature. 
To illustrate the comparison of the methods in terms of relative error, we perform numerical experiments in MATLAB with various matrices. The experiments are performed on a computer with AMD Ryzen 5 4500U having 6 CPUs and 8 GB RAM with OS system Ubuntu 22.04 LTS. Our RQI-PCGTLS-MP algorithm and associated functions are available at \verb|https://github.com/edoktay/rqipcgtlsmp|, which includes
scripts for generating the data and plots in this work.	

Recall that in the $k$-th Rayleigh quotient iteration, both algorithms perform at most $k+1$ PCGTLS iterations in each of lines \ref{algline1} and \ref{algline2} of Algorithm \ref{alg:rqi_mp}. Recall also that we terminate both RQI-PCGTLS and RQI-PCGTLS-MP algorithms when the normalized residual norm increases, i.e., $\psi_{k+1}>\psi_k$, where $\psi_k$ is defined in \eqref{eq:psik}. For a detailed explanation, see Section \ref{sec:rqipcgtlsmp}.

% For Cholesky factorization in RQI-PCGTLS-MP, the low precision Cholesky algorithm from \cite{hp:21} is used with two-sided diagonal scaling with scaling parameter $\theta=0.1$ and diagonal perturbation. 
Each experiment contains a random error matrix and vector. For reproducibility of random arrays for the experiments (except Example 1 in Section \ref{subsec:ex1}), we use the MATLAB command \verb|rng(1)| each time we run the algorithm.

Each plot shows the relative error versus the number of Rayleigh quotient iterations performed. Dashed lines represent the relative error \verb|rerrx| in the approximate solution, i.e., $\|x_k-x\|_2/\|x\|_2$, whereas solid lines represent the relative error \verb|rerrs| in the Rayleigh quotient after $k$ Rayleigh quotient iterations, i.e., $|\sigma_{n+1}^{(k)}-\sigma_{n+1}|/|\sigma_{n+1}|$. The quantity \verb|rerrs| is calculated with $\sigma_{n+1}$ since the Rayleigh quotient iteratively approaches $\sigma_{n+1}^2$. Thus, \verb|rerrs| calculates the convergence of $\sigma_{n+1}^{(k)}$ to $\sigma_{n+1}$. %We say an algorithm converges when the relative error drops below the unit roundoff for the working precision ($\approx 10^{-16}$ in our case). For units roundoff of different precisions, see Table \ref{tab:eps}. 
Red lines show errors in RQI-PCGTLS-MP and blue lines in RQI-PCGTLS. We used the Advanpix toolbox \cite{advanpix} to simulate quadruple precision for the computation of the plotted quantities, which is not included in the performance analysis. 

In our numerical experiments, we used double {precision in MATLAB for the RQI working precision $u$ to be able to compare with RQI-PCGTLS (performed in double precision). For each problem, we evaluate the constraints on $u_q$ given in }{\eqref{constr1}}{ and }{\eqref{constr2}}{, and show that half precision satisfies both constraints for each of the tested problems. To simulate half precision,} we use the \texttt{chop} library and associated functions from \cite{higham2019simulating}, available at \texttt{https://github.com/higham/chop} and \texttt{https://github.com/SrikaraPranesh/LowPrecision\_Simulation}. 

{Although we do not provide a constraint on $u_p$, we found that using a precision less than fp32 in RQI-PCGTLS-MP when $u$ is fp64 gives }\verb|rerrx|{ and }\verb|rerrs|{ much greater than the fp64 unit round-off. Thus, RQI-PCGTLS-MP fails to converge to $x$ and $\sigma^2_{n+1}$. Although single precision for $u_p$ was sufficient for all included examples, we give an example where lower than single precision for $u_p$ does not work well in Section }\ref{sec:bjorck}. {We also note that, although not shown for each experiment, we have also tested the use of Cholesky factorization instead of Householder QR to compute the preconditioner for PCG. In many cases, the use of Cholesky resulted in a greater number of Rayleigh quotient iterations required for convergence for RQI-PCGTLS and/or RQI-PCGTLS-MP. An exception to this is discussed in Section} \ref{sec:vanhuffel}.

\subsection{Example 1: Random matrix}\label{subsec:ex1}
We first test our algorithm on a random matrix $A$ of dimension $(m,n)=(100,60)$ having uniformly distributed random elements. We generated the matrices using the MATLAB command \verb|rand(m,n)|. A vector of ones is used for the right-hand side vector $b$. The error matrix $E$ and vector $e$ are added to the left and right-hand sides, respectively. The final system can be written as $(A+E)x = b+e$, where $E = \epsilon\cdot$\verb|rand(m,n)| and $e = \epsilon\cdot$\verb|rand(m,1)|. For this example, we use the MATLAB command \verb|rng(0)| for reproducibility and $\epsilon$ is chosen as $10^{-6}$.

\begin{figure}[!ht]
		\includegraphics[width=.48\linewidth]{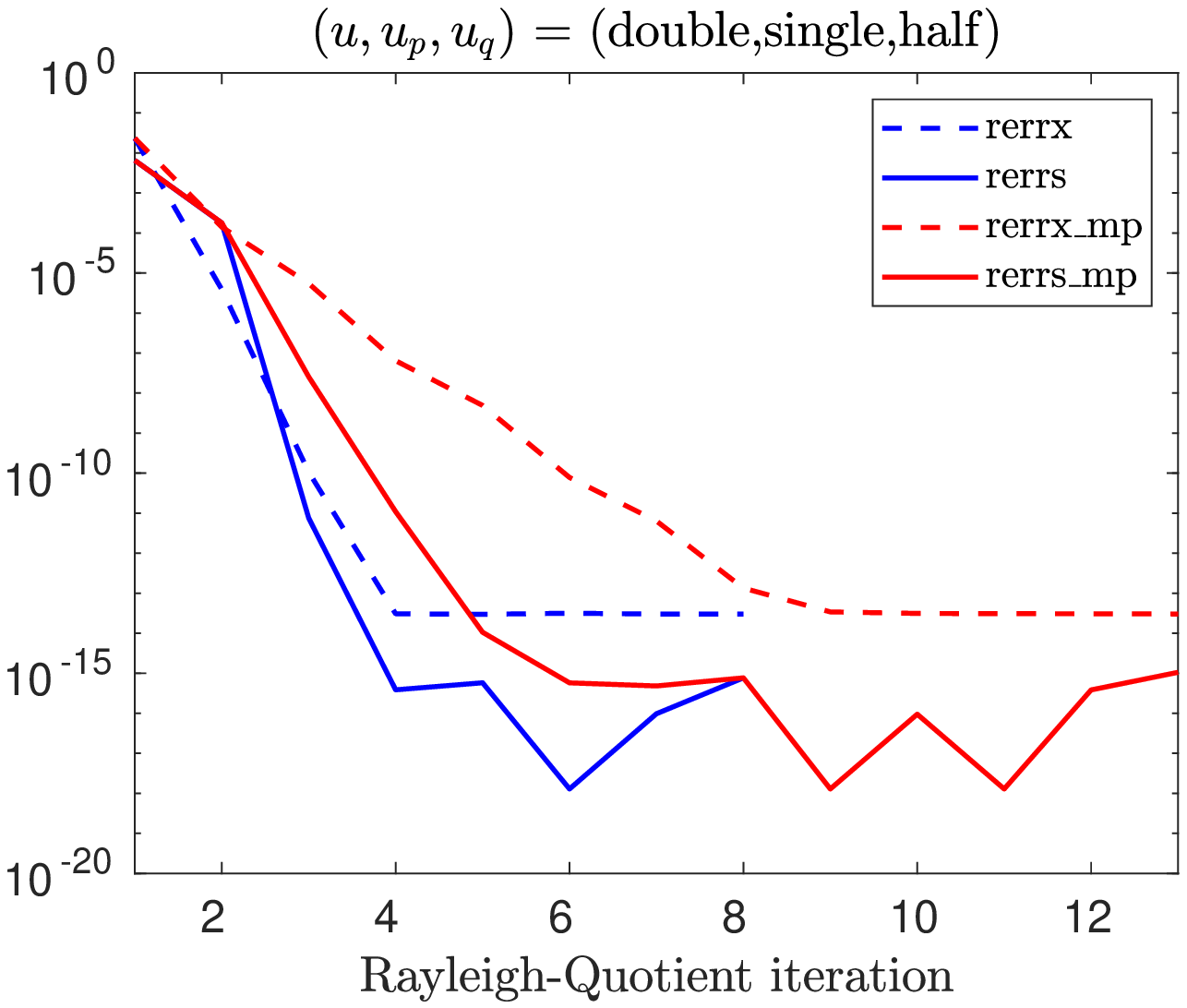}
            \hspace{0.4cm}
		\includegraphics[width=.48\linewidth]{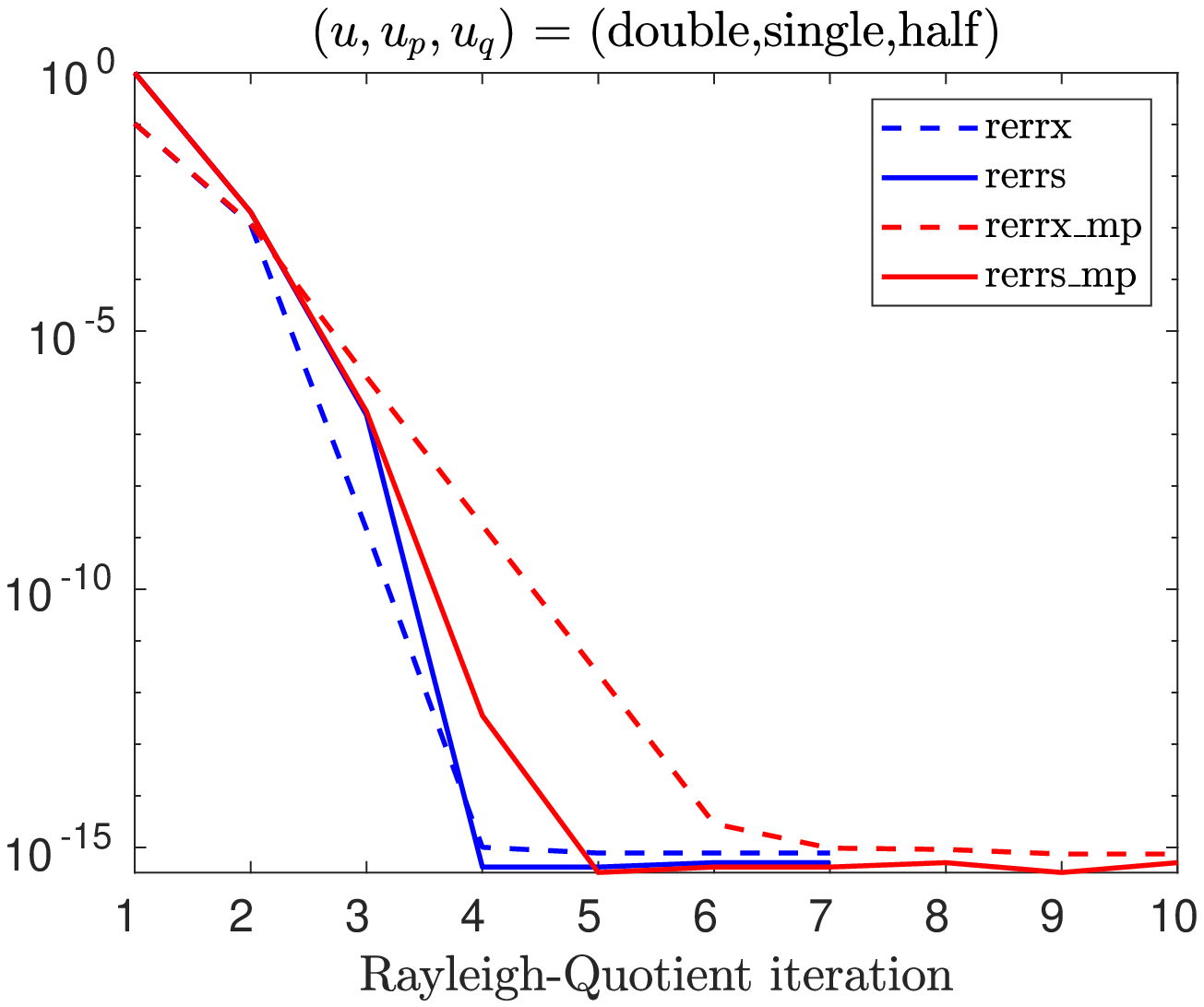}
		\caption{Relative errors versus Rayleigh quotient iterations for the random (left) and $\delta$ (right) matrices.}
		\label{fig:random}
\end{figure}

{For the random matrix, using constraints }{\eqref{constr1}}{ and }{\eqref{constr2}}{, we obtain the bounds $u_q < 2\times 10^{-2}$ and $u_q < 4\times 10^{-3}$, respectively. Thus, half precision is guaranteed to work. }
The left plot in Figure \ref{fig:random} shows the relative error behavior of RQI-PCGTLS and RQI-PCGTLS-MP. We see that RQI-PCGTLS performs 8 Rayleigh quotient iterations, whereas RQI-PCGTLS-MP performs 13 before the stopping criterion is satisfied. 
The mixed precision variant performing more iterations than the uniform precision variant is expected, since we may have a worse preconditioner and thus a worse solution produced by the inner PCG solve for a fixed number of iterations. 

It can be more illustrative to look at the point at which the rate of improvement of the algorithm slows down significantly. For \verb|rerrs|, we see that RQI-PCGTLS improvement starts to slow down at the fourth iteration, and for the RQI-PCGTLS-MP algorithm, it takes six iterations. Importantly, note that the limiting accuracy for both \verb|rerrs| and \verb|rerrx| is similar for both RQI-PCGTLS and RQI-PCGTLS-MP. 

% \begin{figure}[!ht]
% 	\begin{minipage}{60mm}
% 		\includegraphics[width=\textwidth]{rand_total_1.eps}
% 		\caption{Relative errors versus Rayleigh quotient iterations for the random matrix.}
% 		\label{fig:random}
% 	\end{minipage}
% 	\hfill
% 	\begin{minipage}{60mm}
% 		%\vspace*{-0.2cm}
% 		\includegraphics[width=\textwidth]{delta_total_1.eps}
% 		\caption{Relative errors versus Rayleigh quotient iterations for the $\delta$ matrix.}
% 		\label{fig:sigma1}
% 	\end{minipage}
% \end{figure}

\subsection{Example 2: The $\delta$ matrix}
We now test our algorithm on the system $(A+E)x = b+e$ where $A$ is a matrix used in \cite{ds:18}, which we call the $\delta$ matrix. The system $Ax=b$ is constructed as

\begin{eqnarray*}
\begin{bmatrix}
				\delta &0&0&0\\
				0&0&0&0\\
				0&\delta&0&0\\
				0&0&0&0\\
				0&0&0&0\\
				0&0&0&0\\
				0&0&1&0\\
				0&0&0&0\\
				0&0&0&1\\
			\end{bmatrix}\begin{bmatrix}
			x_1\\
			x_2\\
			.\\
			.\\
			.\\
			x_9
		\end{bmatrix}=\begin{bmatrix}
		1\\
		1\\
		.\\
		.\\
		.\\
		1
	\end{bmatrix}.
\end{eqnarray*}
The error matrix $E=\epsilon \bar{E} \cdot A$, and error vector $e=\epsilon \bar{e} \cdot b$ are constructed with $\epsilon=10^{-1}$, $\delta=10^{-2}$, and a uniformly distributed random matrix and vector with entries in the
interval $(-1,1)$, $\bar{E}$ and $\bar{e}$, respectively. For this example, since $\psi_{k+1}\approx\psi_k$ from the tenth iteration onward, we slightly modify the stopping criterion and use $\psi_{k+1}\geq\psi_k$ for this test case.

{For the $\delta$ matrix, using constraints }{\eqref{constr1}}{ and }{\eqref{constr2}}{, we obtain the bounds $u_q < 9\times 10^{-3}$ and $u_q < 1\times 10^{-3}$, respectively. Thus, half precision is guaranteed to work in this case as well. %We also again observed that using Cholesky instead of Householder QR factorization can cause convergence delay in RQI-PCGTLS. The uniform precision algorithm needs to perform one more iteration to reach convergence when used with Cholesky. We note that using Cholesky gives the same result as QR factorization for our RQI-PCGTLS-MP approach. Therefore, we again use Householder QR in both algorithms in this experiment to illustrate error behaviors.
}

The right plot in Figure \ref{fig:random} shows the relative error behavior of RQI-PCGTLS and RQI-PCGTLS-MP for the $\delta$ matrix. Here, RQI-PCGTLS performs 7 Rayleigh quotient iterations before the stopping criterion is met, whereas RQI-PCGTLS-MP performs 10. When we look at the point where the rate of improvement in \verb|rerrs| sufficiently slows down, we see that RQI-PCGTLS-MP requires only one iteration more to have approximately the same error as RQI-PCGTLS. Again, another important observation from the plot is that although the mixed precision variant performs more iterations, both RQI-PCGTLS variants reached a similar limiting accuracy. This again shows that using mixed precision does not cause a loss of accuracy for this problem.

\subsection{Example 3: The Bj\"{o}rck matrix}\label{sec:bjorck}
Next, we use the matrix $A$ that is used to test the convergence properties of RQI-PCGTLS in \cite{bhm:00}, which we call the Bj\"{o}rck matrix. The TLS problem is defined as $(A+E)x = b+e$, where 
\begin{eqnarray*}
A = Y\begin{bmatrix}
				D\\
				0
			\end{bmatrix}Z^T\in \mathbb{R}^{m\times n},
\end{eqnarray*}
$b = Ax$ with $(m,n)=(30,15)$, $x = (1,1/2,\ldots,1/n)$, $D = \text{diag}(1,2^{-1},\ldots,2^{-n+1}),$ and $Y,Z$ are random orthogonal matrices generated using the MATLAB command \verb|RandOrthMat()|. The error matrix $E = \epsilon\cdot$\verb|rand(m,n)| and vector $e = \epsilon\cdot$\verb|rand(m,1)| are generated with $\epsilon = 0.05$.

{For the Bj\"{o}rck matrix, using constraints }{\eqref{constr1}}{ and }{\eqref{constr2}}{, we obtain the bounds $u_q < 4\times 10^{-2}$ and $u_q < 5\times 10^{-3}$, respectively. Thus, half precision is guaranteed to work.}

\begin{figure}[!ht]
		\includegraphics[width=.48\linewidth]{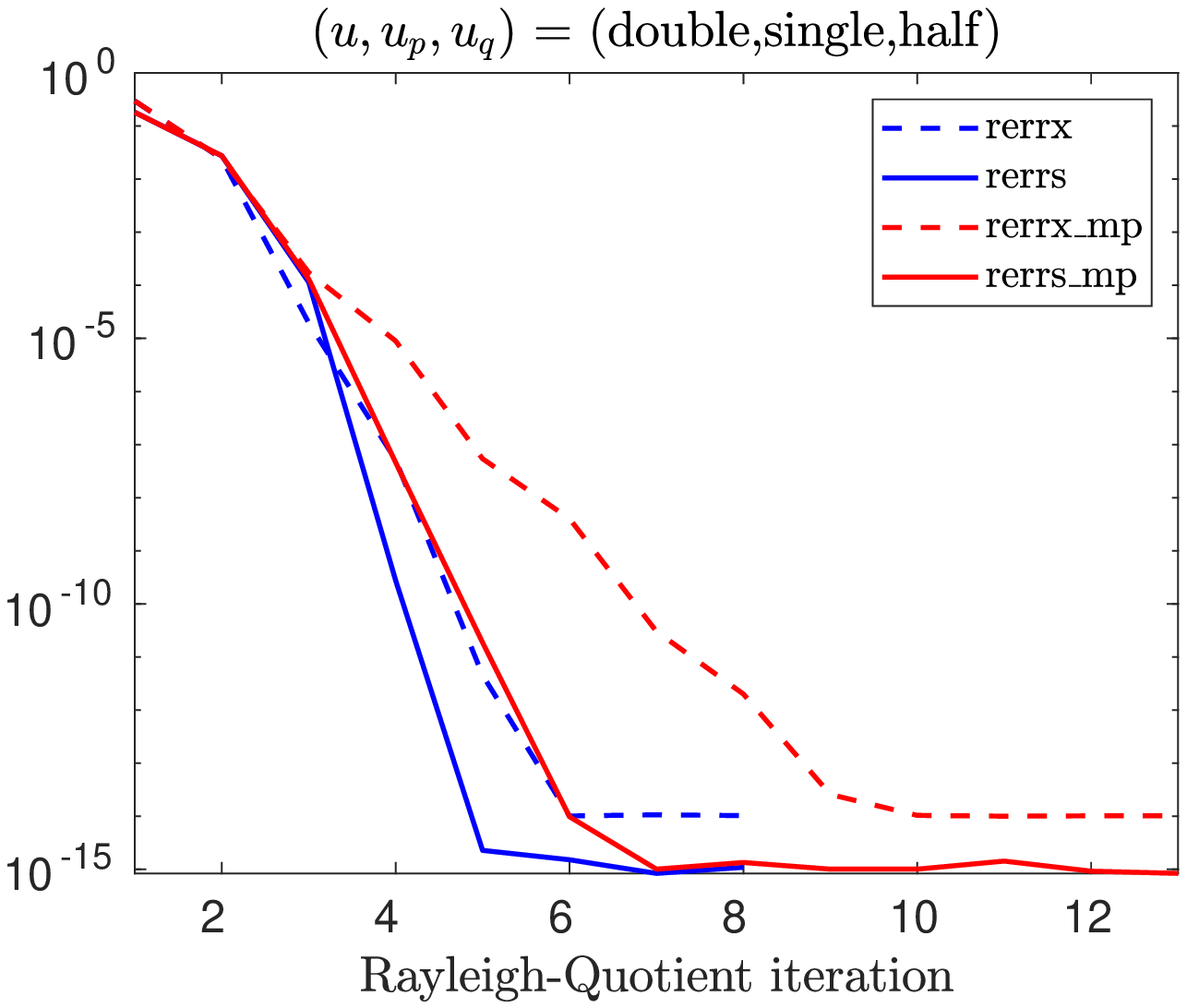}
            \hspace{0.4cm}
		\includegraphics[width=.48\linewidth]{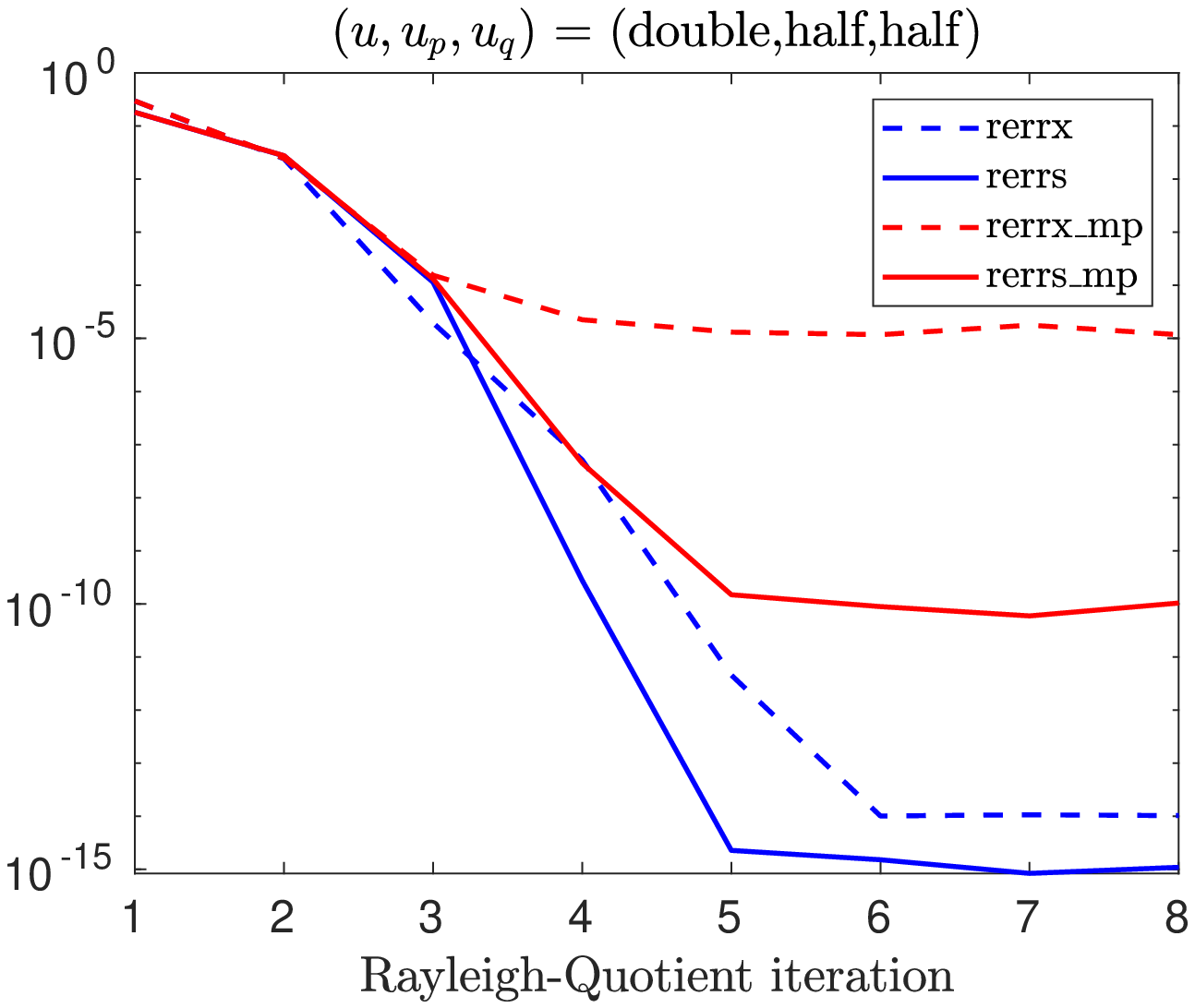}
		\caption{{Relative errors versus Rayleigh quotient iterations for the Bj\"{o}rck matrix using single precision (left) and half precision (right) for $u_p$.}}
		\label{fig:bjorck_half}
\end{figure}

In {the left plot in Figure} \ref{fig:bjorck_half}, {which uses the usual precision settings,} we see similar behavior as in the previous examples. Therefore, we again conclude that while RQI-PCGTLS-MP can result in a convergence delay, both algorithms eventually converge to a similar level. 
{The right plot shows the same experiment, but with $u_p$ set to half instead of single. In this case, while the relative error in the approximate solution given in the dashed line is around $10^{-14}$ for RQI-PCGTLS, it is only around $10^{-4}$ with RQI-PCGTLS-MP.} 

%{We again observed that using Cholesky instead of Householder QR factorization causes convergence delay and additional iterations after convergence until reaching the termination criteria in both algorithms. } 

% \begin{figure}[!ht]
% 	\begin{minipage}{60mm}
% 		\includegraphics[width=\linewidth]{bjorck_total_1.eps}
% 		\caption{Relative errors versus Rayleigh quotient iterations for the Bj\"{o}rck matrix.}
% 		\label{fig:test}
% 	\end{minipage}
% 	\hfill
% 	\begin{minipage}{60mm}
% 		%\vspace*{-0.2cm}
% 		\includegraphics[width=\linewidth]{toeplitz_total_1.eps}
% 		\caption{Relative errors versus Rayleigh quotient iterations for the Toeplitz matrix.}
% 		\label{fig:toeplitz}
% 	\end{minipage}
% \end{figure}

\begin{figure}[!ht]
\centering
		\includegraphics[width=.48\linewidth]{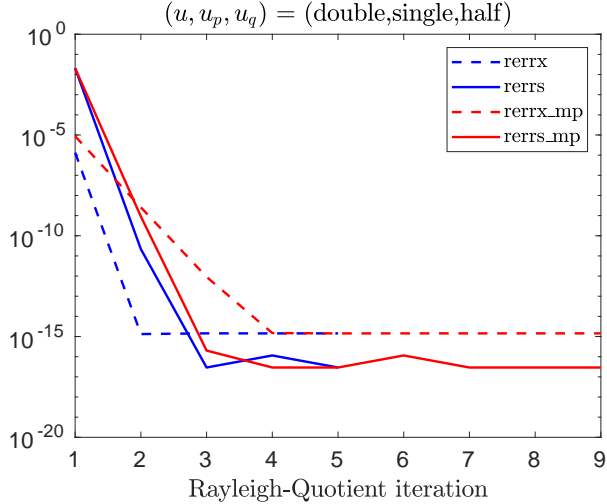}
		\caption{Relative errors versus Rayleigh quotient iterations for {the Toeplitz matrix.}}
		\label{fig:toeplitz}
\end{figure}
\subsection{Example 4: The Toeplitz matrix}
The matrix $A$ is constructed as $A = \bar{T}+E,$ where $E$ is a random Toeplitz matrix and $\bar{T}\in \mathbb{R}^{n\times (n-2\omega)}$ is the lower Toeplitz matrix having elements 
\[t_{i,1} =\begin{cases}
				\frac{1}{\sqrt{2\pi \alpha^2}}e^{[\frac{-(\omega-i+1)^2}{2\alpha^2}]} & i=1,2,\ldots,2\omega+1,\\
				0& \text{otherwise,}
			\end{cases}\]
    in the first row and 
  \[t_{1,j}=\begin{cases}
			t_{1,1} & \text{if }j=1,\\
			0&\text{otherwise}
		\end{cases}\] 
  in the first column. This matrix comes from an application in signal restoration \cite{t:83}. The matrices $\bar{T}$ and $E$ are constructed using the built-in MATLAB function \verb|toeplitz()| with the input of the first column of $\bar{T}$. The right-hand side vector is constructed as $b = \bar{b}+e$, where $\bar{b}$ is a vector of ones, and $e$ is a random vector. For this test, we use $n=100$, $\omega = 2$, and $\alpha=1.25$. We scale $E$ and $e$ with a scaling factor of 0.001.

  {For the Toeplitz matrix, using constraints }{\eqref{constr1}}{ and }{\eqref{constr2}}{, we obtain the bounds $u_q \lesssim 1$ and $u_q \lesssim 10^{-2}$, respectively. Thus, half precision is guaranteed to work here. 
  %We also observed that we need two more iterations to reach convergence when Cholesky was used in the RQI-PCGTLS algorithm, while RQI-PCGTLS-MP was not affected by the change of the factorization method. Therefore, we again use Householder QR in both algorithms in this experiment.
  }

  Even though {Figure} \ref{fig:toeplitz} shows that \verb|rerrs| oscillates after the third iteration for both uniform and mixed precision variants, the oscillations are negligible since they happened after the error reaching the level indicated by the working precision, $u\approx10^{-16}$. We also see that, as in previous examples, both variants converge to a similar level although mixed precision may result in slightly delayed convergence.
  
\subsection{Example 5: The Van Huffel matrix}
\label{sec:vanhuffel}

Finally, we test the matrix with $n=100$ from \cite{vv:91}, which is used to illustrate the difference between LS and TLS solutions. We define the TLS problem by $(A+E)x = b+e$, where $Ax = b$ is given as
\[ \underbrace{\begin{bmatrix}
						n-1&-1&\cdots&-1\\
						-1&n-1&\cdots&-1\\
						.&&&\\
						.&&&\\
						.&&&\\
						-1&-1&\cdots&n-1\\
						-1&-1&\cdots&-1\\
						-1&-1&\cdots&-1\\
				\end{bmatrix}}_{n\times n-2}\begin{bmatrix}
					x_1\\
					.\\
					.\\
					.\\
					x_{n-3}\\
					x_{n-2}
				\end{bmatrix}\approx \begin{bmatrix}
					-1\\
					-1\\
					.\\
					.\\
					.\\
					-1\\
					n-1\\
					-1
				\end{bmatrix},\]
and $E = \epsilon\cdot$\verb|rand(n,n-2)|, $e = \epsilon\cdot$\verb|rand(n,1)| are the error matrix and vector generated with $\epsilon = 10^{-6}$, respectively.

{For the Van Huffel matrix, using constraints }{\eqref{constr1}}{ and }{\eqref{constr2}}{, we obtain the bounds $u_q < 2\times 10^{-1}$ and $u_q < 5\times 10^{-3}$, respectively. Thus, half precision is guaranteed to work.}

\begin{figure}[!ht]
		\includegraphics[width=.48\linewidth]{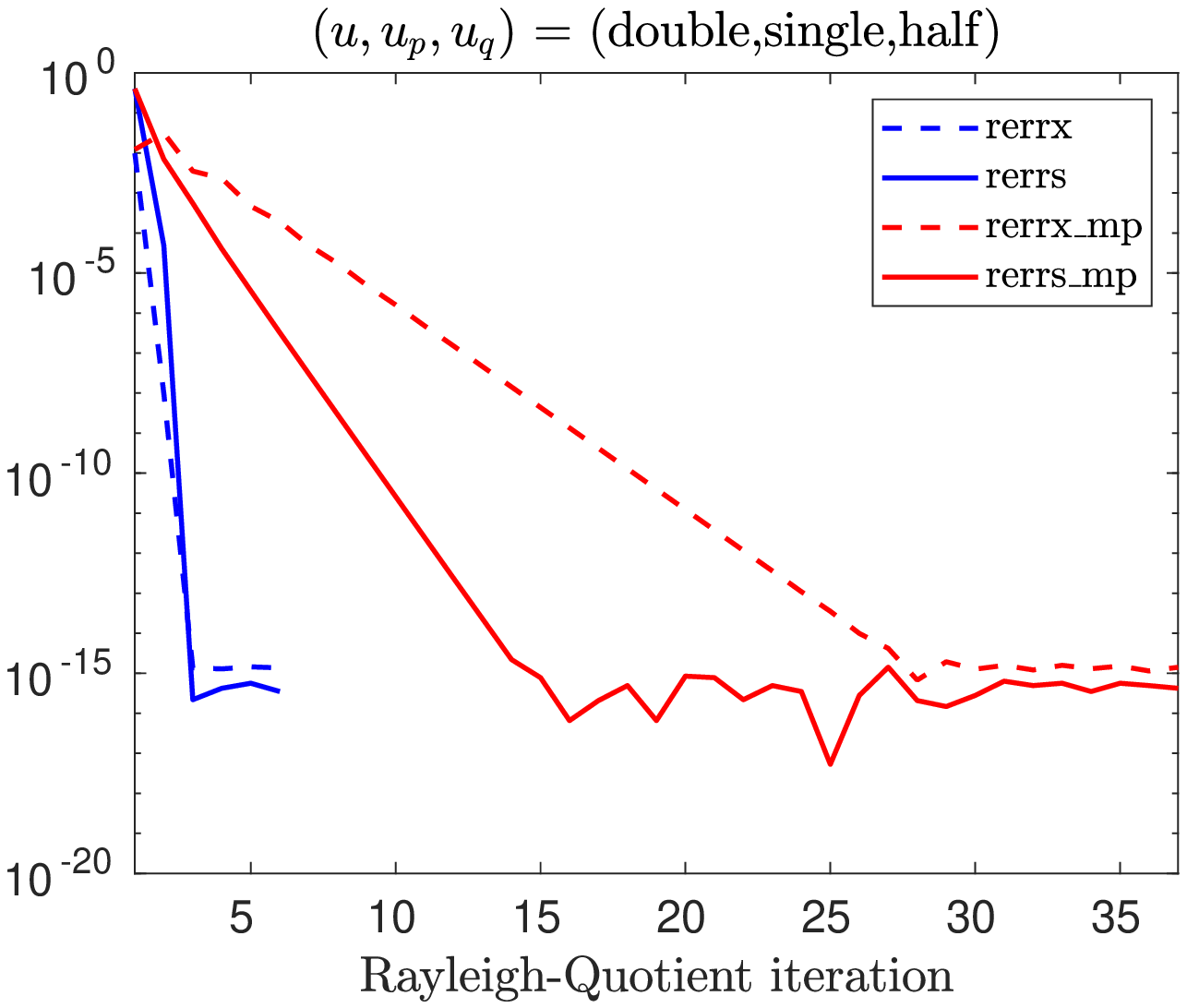}
            \hspace{0.4cm}
		\includegraphics[width=.48\linewidth]{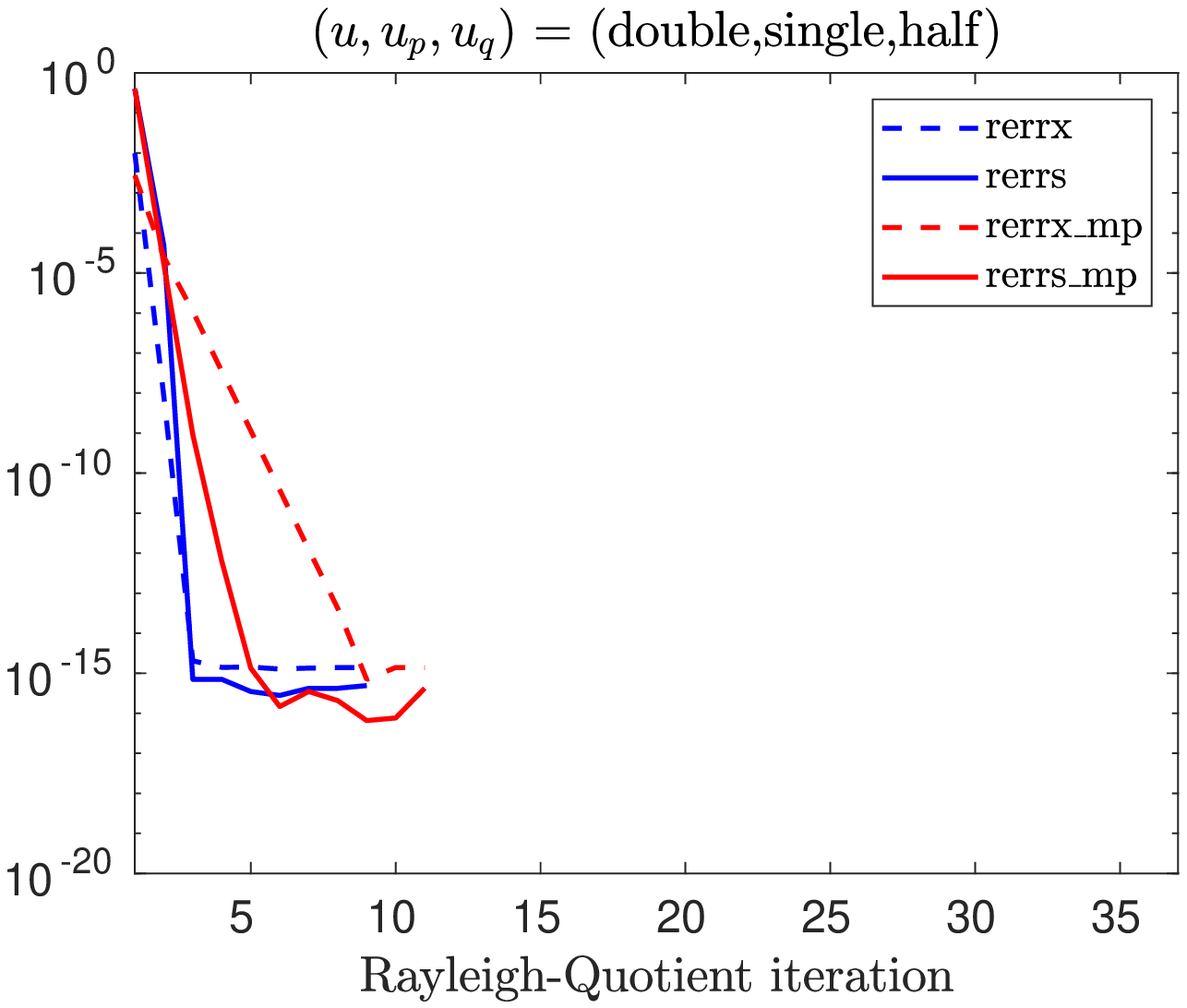}
		\caption{Relative errors versus Rayleigh quotient iterations for the Van Huffel matrix with QR (left) and Cholesky factorization with two-sided diagonal scaling (right).}
		\label{fig:vanhuffel}
\end{figure}

We see exceptional behavior on the RQI-PCGTLS-MP algorithm for this example in the left plot of Figure \ref{fig:vanhuffel}. Here, we observe that using mixed-precision causes the algorithm to perform over $5\times$ more Rayleigh quotient iterations than RQI-PCGTLS if the QR factorization of $A$ is used as a preconditioner in PCGTLS (although at the end, they again both converge to a similar level). For the comparison of RQI-PCGTLS and our mixed precision variant with this matrix, we also tried using the Cholesky factorization as a  preconditioner within both variants. For low-precision Cholesky factorization with two-sided diagonal scaling, we use the \verb|cgir3.m| function available in \texttt{https://github.com/SrikaraPranesh/Multi\_precision\_NLA\_kernels}. In our experiments, we see that for this example, using the Cholesky factor of $A^TA$ is ultimately better than using the QR factorization within RQI-PCGTLS-MP; from the plot on the right, the mixed precision variant with the Cholesky factorization performs only 3 and 6 more iterations to converge to a similar level as RQI-PCGTLS for \verb|rerrs| and \verb|rerrx|, respectively. Thus it is clear that for some matrices, Cholesky factorization makes a better preconditioner than QR factorization within RQI-PCGTLS-MP. A full explanation for this behavior is left as future work.

% \begin{figure}[!ht]
% 	\begin{minipage}{60mm}
% 		\includegraphics[width=\linewidth]{vanhuffel_total_1.eps}
% 		\caption{Relative errors versus Rayleigh quotient iterations for the Vanhuffel matrix with QR factorization.}
% 		\label{fig:vanhuffel}
% 	\end{minipage}
% 	\hfill
% 	\begin{minipage}{60mm}
% 		%\vspace*{-0.2cm}
% 		\includegraphics[width=1.1\linewidth]{vanhuffel_total_chol_1.eps}
% 		\caption{Relative errors versus Rayleigh quotient iterations for the Vanhuffel matrix with Cholesky factorization with two-sided diagonal scaling.}
% 		\label{fig:vanhuffel_chol}
% 	\end{minipage}
% \end{figure}

\section{Conclusion}\label{sec:conclusion}

In this paper, we present a mixed precision variant of RQI-PCGTLS called RQI-PCGTLS-MP. Our algorithm allows the use of up to three different precisions in different parts of the algorithm. To allow the use of lower precision for more ill-conditioned systems, we advocate for the use of the Householder QR factorization of $A$ instead of the Cholesky factorization of $A^TA$ to construct a low-precision preconditioner for PCGTLS. We derive constraints on the precision with which the QR factorization can be computed. Interestingly, in contrast with the standard least squares case, the factorization precision we can use depends not only on $A$ but on the right-hand side $b$ in the total least squares setting. 

The numerical experiments we performed show that using mixed precision increases the number of Rayleigh quotient iterations performed; however, since we are using half precision in the QR factorization and single precision in PCGTLS, the computation cost is reduced per iteration. Thus, from the performance model we constructed, we see that we can compensate for the additional Rayleigh quotient iterations up to $4 \times$ depending on the size of the matrix and how many Rayleigh quotient iterations are performed. Furthermore, our experiments indicate that the use of low precision does not impact the final attainable accuracy of the TLS solution.

Our future work involves a complete mathematical analysis, which we suspect can be done by adapting the results on inexact Newton in \cite{hm:22,t:01,k:22} and inexact RQI in \cite{fs:07,se:02}. With the complete analysis, we will be able to prove our observations on the attainable accuracy and provide a theoretical explanation of different convergence delays on different systems. Furthermore, it would be of great interest to compare the practical run-time of both algorithms on {large-scale} matrices in a high-performance setting with hardware that supports half precision computation. 
	
	\clearpage
	\bibliographystyle{siamplain}
	\bibliography{references}

\end{document}